\documentclass[preprint,number, 12pt]{elsarticle}

\usepackage{graphicx}
\usepackage{amssymb,amsmath,amsthm}

\usepackage{trackchanges}
\usepackage{soul}
\usepackage{caption}
\usepackage{subcaption}
\usepackage{trackchanges}
\usepackage{url}
\usepackage{epstopdf}
\usepackage{cleveref}
\usepackage{soul,color}

\newcommand{\hl}[1]{#1}
\newcommand{\hlblue}[1]{#1}

\newcommand{\crefrangeconjunction}{--}
\newtheorem{theorem}{Theorem}

\newtheorem{lemma}{Lemma}
\newtheorem{proposition}{Proposition}

\theoremstyle{definition}

\newcommand{\mA}{{\mathcal A}}

\newcommand{\mI}{{\mathcal I}}
\newcommand{\mJ}{{\mathcal J}}

\newcommand{\mL}{{\mathcal L}}

\newcommand{\mS}{{\mathcal S}}

\newcommand{\mY}{{\mathcal Y}}

\newcommand{\bR}{{\mathbb R}}
\newcommand{\bZ}{{\mathbb Z}}

\newcommand{\bE}{{\mathbb E}}
\newcommand{\bI}{{\mathbb I}}
\begin{document}

\begin{frontmatter}

\title{Decentralized Signal Control for Urban Road Networks
\author[add2]{Tung~Le}\corref{cor1}
\ead{tmle@swin.edu.au}
\author[add1]{P\'{e}ter~Kov\'{a}cs}
\ead{P.Kovacs@uva.nl}
\author[add1]{Neil~Walton}
\ead{n.s.walton@uva.nl}
\author[add2]{Hai~L.~Vu}
\ead{hvu@swin.edu.au}
\author[add2]{Lachlan~L.~H.~Andrew}
\ead{l.andrew@ieee.org}
\author[add3]{Serge~S.~P.~Hoogendoorn}
\ead{s.p.hoogendoorn@tudelft.nl }
}
\cortext[cor1]{Corresponding author}
\address[add2]{Faculty of ICT, Swinburne University of Technology, Australia}
\address[add1]{Korteweg-de Vries Institute for Mathematics, Universiteit van Amsterdam, The Netherlands}
\address[add3]{Delft University of Technology, Delft, The Netherlands.}

\begin{abstract}
We propose in this paper a decentralized traffic signal control
policy for urban road networks. Our policy is an adaptation of a
so-called BackPressure scheme which has been widely recognized in
data network as an optimal throughput control policy. We have
formally proved that our proposed BackPressure scheme, with fixed
cycle time and \hl{cyclic phases}, stabilizes the network for any
feasible traffic demands. Simulation has been conducted to compare
our BackPressure policy against other existing distributed control
policies in various traffic and network scenarios. Numerical results
suggest that the proposed policy can surpass other policies both in
term\hlblue{s} of network throughput and congestion.

\end{abstract}


\end{frontmatter}

\section{Introduction}




Traffic congestion is a major problem in modern societies due to
increasing population and economic activity. This motivates the need
for better utilizing the existing infrastructures and for
efficiently controlling the traffic flow in order to minimize the
impact of congestion.

One of the key tools for influencing the efficiency of traffic flow
in urban networks is traffic signal control that enables conflicting
traffic to flow through intersections via the timing of green/red
light cycles. It has long been recognized that the challenge is to
find optimal cycle timing over many intersections so as to reduce
the overall congestion and to increase the throughput through the
network.

There has been much work in the past both on designing and
optimizing isolated or coordinated signals that reactively resolve
congestion in the urban networks. Broadly, there are two types of
control that have been used for signal control: static and
vehicle-actuated controls; see \cite{Hamilton2013}. Static control
(sometimes called ``fixed time plan'') involves the optimization of
the cycle time, the offset between nearby intersections,
\textit{and} the split of green times in different directions within
a cycle. This can be optimized in isolation or in \hlblue{a coordinated}
manner, for instance to create a so-called green wave where vehicles
always arrive at intersections during the green cycle time, e.g.
\cite{Webster1958,Gartner1975a,Gartner1975b,Kraft2009}. In contrast,
vehicle-actuated controls use online measurements from on-road
detectors (e.g., inductive loops) to optimize signal timings on a
cycle-to-cycle basis in real time. Some examples of commonly used
implementations are: SCOOT \cite{hunt1981scoot}; UTOPIA
\cite{mauro1990utopia}; and the hierarchical scheme RHODES
\cite{Mirchandani2001}. Combinations of both the fixed time plan and
vehicle-actuated control also exist; one widely used example is
SCATS \cite{lowrie1982sydney}.

Given a choice of the control scheme, various approaches to optimize
the signal plans have been proposed. Examples include Mixed-Integer
Linear Programming problems, see
\cite{Gartner1975a,Gartner1975b,Dujardin2011}; Linear Complementary
Problem, see \cite{DeSchutter1999}; rolling horizon optimization
using dynamic programming, see \cite{gartner1983opac, Farges1983,
Mirchandani2001}, or its combination with online learning algorithms
, see \cite{Cai2009}; store-and-forward models based on Model
Predictive Control (MPC) optimization \cite{Aboudolas2009,
Aboudolas2010, Tettamanti2008, Tettamanti2010a, Tettamanti2010b,
Le2013}, or MPC optimization with non-linear
prediction~\cite{Lin2011}.  Many of these approaches formulate the
problem in a way that is centralized and thus \hlblue{are} inherently not
scalable. \hl{While the state of the art is the use of centralized
techniques, improved scalability may be obtained using decentralized
approaches. In this paper, we focus exclusively on decentralized
schemes. Although such schemes are in their infancy, this research
is one step along the path of improving such schemes to offer
performance comparable with centralized schemes while retaining
their scalability.}

A scalable \hl{distributed approach} is to solve a set of loosely coupled
optimizations, one for each intersection, with coupling provided by
traffic conditions.  Two natural approaches are to control the
traffic lights based on either (a) the expected number of vehicles
to enter the intersection during the next cycle, or (b) the
\emph{difference} in traffic load on the road leading into the
intersection and those leading out. These approaches are now
deployable in practice thanks to emerging technologies, such as
cameras and wireless communication enabling better access to
real-time traffic data.

Notable among the first class is the work of \cite{Smith1980}, a
so-called $P_0$ policy and its variants \cite{MUSIC2000,Smith2011},
followed by the work of L\"ammer and Helbing \cite{Lammer2008} where
the switching cost between phases is taken into account. In this
approach, each intersection estimates the amount of traffic that
will arrive during the next complete cycle, and sets the split time
such that each phase gets a time proportional to the number of cars
expected to arrive on roads which have a green light during that phase.
The lack of central control raises the possibility that
intersections may interact in unexpected ways to cause instability.
To limit this, a stabilization mechanism was proposed
in~\cite{Lammer2010}. However, beyond heuristic arguments, there
remains no formal proof of stability of this approach.

Approach (b) including work by \cite{varaiya2013max},
\cite{Wongpiromsarn2012} and \cite{Zhang2012} was inspired by
research developed for packet scheduling in wireless networks: a
so-called max weight or back pressure (refer to as BackPressure in
this paper) algorithm \cite{tassiulas1992stability, McKeown1999}.
Like approach (a), BackPressure does not require any \textit{a
priori} knowledge of the traffic demand, but it has the added
benefit of provable stability. To make that more precise, define a
traffic load to a network as ``feasible'' if there exist splits at
each intersection such that the queues do not build up indefinitely.
Under certain simplifying assumptions, it can be shown that the
queues under BackPressure do not build up indefinitely for any
feasible traffic load. This will be made more formal in
Section~\ref{sec:math}. In wireless networks, BackPressure can be
computationally prohibitive, but in road networks
\cite{Wongpiromsarn2012, varaiya2013max, Zhang2012}, it admits a
simple distributed implementation, just like approach~(a).

\hl{It is worth noting that all the above mentioned policies
\cite{Smith1980, Smith2011, Lammer2008, varaiya2013max,
Wongpiromsarn2012, Zhang2012} \hlblue{make decisions periodically bases on
the evaluations of traffic over a fixed time interval.  These are called \emph{fixed cycle} policies.} For example, the BackPressure policy
\cite{Wongpiromsarn2012} determines the phase to be activated at the
beginning of each fixed time slot, while the policy
\cite{Lammer2010} decides whether to keep serving the current flow
or switch to other flow at a regular time interval which can be
arbitrary small.}

Given the possibility of a stability guarantee by the BackPressure
scheme, our objective in this work is to fully adapt it to the
traffic control scenarios. To this end, we propose in this paper a
new signal control strategy that addresses two weaknesses in the
prior application of BackPressure to road networks
\cite{Wongpiromsarn2012, varaiya2013max, Zhang2012} while retain and
prove the important stability property of the BackPressure-based
algorithms.

The first \hl{weakness to be addressed is that phases can form an
erratic, unpredictable order in the previously proposed BackPressure
policy.} This is acceptable in the context of communications systems
but for urban road traffic this is undesirable since erratic
ordering of phases brings frustration to drivers and potentially
causes confusion leading to dangerous actions. Moveover, if one
inbound road is particularly backlogged, then it is possible that
other roads are ``starved'' by being assigned a red light for an
extended period. To rectify this, we modify BackPressure to a
\hlblue{``cyclic phase''} \hl{policy where a policy is said to be cyclic
phase policy if it allocates strictly positive service time to all
phases in each control decision, and thus, it is possible to arrange
the phases into a fixed ordered sequence.}

The second \hl{weakness} that we address is that prior applications have
required each intersection to know the ``turning fractions'', that
is, the fraction of traffic from each inbound road that will turn
into each possible outbound road.  We prove that the stability
results still apply when these turning fractions are estimated using
even very simple measurements; specifically, any unbiased estimator
of the turning fractions suffices. Such stability proofs apply for a
general network model but under idealized assumptions. Nonetheless
these form an important step towards the application of BackPressure
to real networks.

To test the practicality of the theoretical refinements described
above, we also present the numerical comparison of the proposed
BackPressure algorithm with the approach of \cite{Smith1980} or
\cite{Lammer2008,Lammer2010} without switching cost. \hl{The}
results suggest that our cyclic phase BackPressure policy tends to
outperform other \hlblue{distributed} polices both in terms of throughput and
congestion. Although the performance of each policy varies widely
depending on the parameter setting such as cycle length or decision
frequency, under the optimal setting, the BackPressure with cyclic
phase and without cyclic phase
 have better throughput in compare with the other policies.

The rest of this paper is organized as follows. We first present the
notations and queue dynamics model before describing our proposed
cyclic phase BackPressure policy in Section~\ref{sec:policy}. The
main results for stability of our policy are then provided in
Section~\ref{sec:math}. These results in a certain sense mean, that
we can interpret our policy as \hlblue{stabilizing the system for the largest} possible set of arrival rates leading to sufficient throughput even in
congested network. For readability, however, most of the
mathematical details and derivations are listed in the Appendices of
the paper. Section~\ref{sec:num} presents the simulation results and
numerical comparison of our scheme with other existing policies
where we demonstrate the benefits of the proposed cyclic phase
BackPressure signal control strategy. Finally,
Section~\ref{sec:concl} concludes the paper and discusses future
work.

\section{\hl{Cyclic Phase BackPressure} Traffic Signal Control}\label{sec:policy}

\subsection{Notation and Network description}\label{sec:sharedNotation}

Consider a network of traffic intersections. This road network
consists of a number of \emph{junctions}, indexed by $\mJ$. Each
junction $j\in\mJ$ consists of a number of \emph{in-roads}, $\mI_j$.
Note that the $\mI_j$ are mutually disjoint, and denote $\mI =
\cup_{j\in\mJ} \mI_j$. A road with multiple lanes having different
turning options (such as a left-turn only lane) is modeled as
multiple in-roads, thus an in-road may model one or more lanes of
traffic flow. Whether these traffic flows are conflicting or not is
not considered in this setting. We use the inclusion $i\in j$ to
indicate that in-road $i$ is part of junction $j$, and we let $j(i)$
notate the junction used by in-road $i$.

Each junction may serve different combinations of in-roads
simultaneously. We call a combination of in-roads served
simultaneously a \emph{service phase}. A service phase for junction
$j$ is represented by a vector $\sigma=(\sigma_i : i\in j)$ where
$\sigma_i$ denotes the rate at which cars can be served from in-road
$i$ at junction $j$ during phase $\sigma$.  In particular,
$\sigma_i>0$ if in-road $i$ has a green light during phase $\sigma$, or
$\sigma_i=0$ otherwise.

Let $\mS_j$ denote the set of phases at junction $j$. We will let
$\mL$ denote the set of links of the road traffic network.
\hlblue{Each link represents a road connecting the junctions of the
urban road network}. Here we write $ii'\in\mL$ if it is possible for
cars served at in-road $i$ junction $j$ to next join in-road $i'$
junction $j'$.

In the rest of this section we impose the additional constraint that
all junctions have a common \emph{cycle length} \hl{$T$}, the time devoted to
serving cars from the different in-roads at the junction. Thus we can model time as discrete and consider a slotted time model where $t=0,1,2,...$ denotes the number of the cycle about to
be initiated.
\hl{Control decisions in our policy are then made at the beginning of each time slot~(so it is a fixed cycle policy which is similar to the policies in \cite{Lammer2010, varaiya2013max, Smith1980}).
We also assume that all junctions have the same loss of service $L$ due to idle times during switches and setups. Thus} at each time step, the system decides at each junction $j$
how much time within the next interval to spend serving each phase $\sigma\in\mS_j$ with the constraint that each service phase must be enacted
for some non-zero length of time \hl{and that the sum of the allocated times must not be greater than $T-L$}. We
assume that a car served at one junction in one time interval presents at an in-road of the next junction in its route at the next time interval.

\subsection{Queue Dynamics Model}
Let the queue length $Q_i(t)$ denote the number of cars at in-road
$i\in\mI$ at the beginning of the $t$th traffic cycle, \hl{and denote the vector of queue lengths by $Q(t)=(Q_i(t): {i\in\mI}) $. The decisions in the policy will be based on the measured queue lengths $\hat{Q}(t)$, which might differ from the actual value of $Q(t)$ by an error term as described in the following equation.
\begin{equation}\label{hat_error}
\hat{Q}_i(t)=Q_i(t)+\delta_i(t),
\end{equation}
where the error term $\delta_i(t)$ is bounded and independent of $Q_i(t)$ or the terms at other in-roads. We denote the vector of the error terms by $\delta(t)=(\delta_i(t) :{i\in\mI})$.}

Let $P^j_{\sigma}(t)$ denote the proportion of the traffic cycle at
junction $j$ which is devoted to service phase $\sigma$. For any
policy and for all $j\in\mJ$, we require\hl{
\begin{equation}
 \sum_{\sigma\in\mS_j} P^j_{\sigma}(t)  = 1-\frac{L}{T}\qquad\text{and}\qquad P^j_\sigma(t)>0.
\end{equation}}
\hl{Recall that $\sigma_i$ gives the rate at which cars can be served}
at in-road $i$ if the entire traffic cycle \hlblue{were devoted to phase $\sigma$,}
and $P^j_{\sigma}(t)$ gives the proportion of the traffic cycle
devoted to service phase $\sigma$. \hl{Their product}, $\sigma_i
P^j_\sigma(t)$ gives the expected number of cars to leave in-road
$i$ under service phase $\sigma$, provided the in-road is not
emptied. Accordingly if we let the \hl{random variable} $S_i(t)$ be the potential number
of cars served from in-road $i$ at junction $j$ in traffic cycle
\hl{$t$, the} mean of $S_i(t)$ must satisfy
\begin{equation}\label{Scond}
\bE [ S_i(t) | Q(t),\delta(t) ] =\sum_{\sigma \in \mS_j} \sigma_i P^j_\sigma(t),
\end{equation}
\hl{where we note that the proportions  $P^j_{\sigma}(t)$ are
allocated according to the decision in the policy based on
$\hat{Q}(t)$, thus the dependency on $\delta(t)$.} The \hl{random}
variable $S_i(t)$ only gives the number of cars served if the
junction does not empty. Thus, it may be possible for $S_i(t)$ to be
greater than the queue size $Q_i(t)$. In this case, $Q_i(t)$ will be
the number of cars served. In other words, the number of cars
actually served at junction $i$ is
\begin{equation}
S_i(t)\wedge Q_i(t)
\end{equation}
where $x\wedge y= \min \{ x,y\}$.

Further, when traffic is served it will move to neighbouring
junctions. For $ii'\in\mL$, we let $p_{ii'}(t)$ denote the
proportion of cars served at in-road $i$ that subsequently join
in-road $i'$ at time $t$. We assume that cars within an in-road are
homogeneous in the sense that each car at the junction has the same
likelihood of joining each subsequent junction. We denote the
expectation of $p_{ii'}(t)$ by $\bar{p}_{ii'}$. We further assume
that this likelihood is constant \hlblue{(i.e. time independent)}
and will not be altered by the queue lengths observed by cars within
the network. Thus $[S_i(t)\wedge Q_i(t)] p_{ii'}(t)$ is the number
of cars that leave inroad $i$ and, next, join inroad $i'$ provided
the in-road does not empty.

We let $A(t)=(A_{i}(t):i\in\mI)\in\bZ_+^\mI$ denote be the number of
external arrivals at in-road $i$ at time $t$. \hl{The expected
number of arrivals or \emph{arrival rate} into each in-road at time
$t$ is defined as $\bar{a}_i(t):=\bE[A_i(t)]$. Notice by allowing
$\bar{a}_i(t)$ to vary as a function of time, we can model varying
traffic demands which undoubtedly can change over the course of a
day. In cases where we choose arrival rates to be static and
unchanging with time, then we will simply denote these arrival rates
by $a_i$.}

Given a service policy $\{ P(t) \}_{t=0}^\infty$, we can define the
dynamics of our queueing model. In particular, we define for
in-road $i$ of junction $j$
\begin{equation}\label{Qprocess}
Q_i(t+1) = Q_i(t) - S_i(t)\wedge Q_i(t)   + A_i(t) + \sum_{i' : i'i\in\mL }[S_{i'}(t)\wedge Q_{i\rq{}}(t)]   p_{i'i}(t).
\end{equation}
Here we assume that cars first depart within a traffic cycle and
then subsequently cars arrive from other in-roads.

\subsection{\hl{Cyclic Phase} BackPressure Control Policy}\label{sec:fixed.BP}
Now we are ready to give our proposed policy as follows

\begin{enumerate}

\item
At the beginning of each traffic cycle, \hl{form an estimate of the
actual turning fractions $p_{ii'}(t)$ with the unbiased estimator
$\bar{q}_{ii'}(t)$.\footnote{\ref{sec:apx.est} gives an expanded
explanation of the estimation method and proposes a form of
$\bar{q}_{ii'}(t)$.}}
\item
For each junction $j\in\mJ$, calculate the weight associated with
each service phase at the junction as a function of the measured
queue sizes $\hat{Q}(t)$ and the above defined estimated
\hl{turning} probabilities \hl{
 \begin{equation}\label{wdef}
w_{\sigma}(\hat{Q}(t)) =  \sum_{i\in j} \sigma_i \left( \hat{Q}_i(t) - \sum_{i': ii'\in\mL} \bar{q}_{ii'}(t)  \hat{Q}_{i'}(t)\right).
\end{equation}}
\item Given these weights, assign the following proportion of the common cycle length to each phase $\sigma$ in $\mS_j$ within the next service cycle,
\begin{equation}\label{policydef}
P_\sigma^j(t)= \frac{\exp\left\{\eta w_{\sigma}(\hat{Q}(t))\right\} }{\sum_{\pi\in \mS_j}  \exp\left\{\eta w_{\pi}(\hat{Q}(t))\right\} },
\end{equation}
for $j\in\mJ$ and where $\eta>0$ is a parameter of the model.\footnote{Given this definition the weights do not need to be strictly positive.}
\end{enumerate}

\hl{The weights defined in \eqref{wdef} are used in the BackPressure policy as
given by \cite{tassiulas1992stability}. They can be viewed as a ``pressure'' a queue places
on downstream queues, which is given by the weighted \hlblue{mean of the} difference\hlblue{s} of
the queue sizes. The larger the weight associate with \hlblue{a phase},
the more important it is \hlblue{to serve the in-roads with green lights during that phase}. Then those weights
are used to calculate the portion of the traffic cycle for each
phase according to \eqref{policydef}.
\hlblue{The distribution \eqref{policydef} gives each phase positive service, with more service given to the higher weight phases. As $\eta \rightarrow 0$, the service allocation tends to uniform, and as $\eta \rightarrow \infty$, the fraction of service given to the highest weight phase(s) tends to $1$.}}

\hl{Notice that in contrast to BackPressure policies which always serve
the phase associated with the highest weight, the proposed policy
ensures that each phase (and subsequently each in-road) receives
non-zero service in each and every cycle. Thus this ensures a cyclic
phase policy, while maintaining the property that higher weights
result in higher proportions of  allocated green time. Note the
policy can be implemented in a decentralized way, after each
junction has communicated queue sizes with its upstream in-roads,
the phases can be calculated. This decentralization has numerous
advantages: it is computationally inexpensive, it does not require
centralized aggregation of information and thus is easier to
implement, and it increases the road networks robustness to
failures.}

\hl{The cyclic phase feature which we introduce to the BackPressure
polices is important from the users' point of view for various
reasons. Firstly, the drivers usually expect an ordered phase
sequence and anticipate traffic signal changes in advance. Secondly,
the waiting time to receive service for any in-road is bounded in
our policy while it could be arbitrary large for some in-roads in
previous \hlblue{state-of-the-art distributed policies, such as BackPressure}. It is particularly important when
considering that pedestrian phases might also be initiated in
parallel with the service phases for vehicles.}

\hl{From an implementation point of view the policy is desirable since
it does not} require knowledge regarding the destination of each car
within the road network, nor does it assume that the proportion of
cars moving between links is known in advance. The policy estimates
turning fractions and measures queue sizes in an on-line manner, and
uses this adaptive estimates and the measurement results to inform
the policy decision.

\section{Mathematical Results - Stability of Cyclic Phase BackPressure Control Policy}\label{sec:math}

\subsection{Stability Region and Queueing Stability}
We define the stability region \hl{$\mA$} of the network to be the set \hl{of arrival rate vectors}
$a=(a_i: i\in\mI)\geq 0$, for which
there exists a positive vector $\rho=(\rho^j_{\sigma}:
\sigma\in\mS_j, j\in\mJ)$, namely the green time proportion devoted
to the service phases in a cycle, and a positive vector
$s=(s_i:i\in\mI)$, namely the departure rates, satisfying the
constraints\hl{
\begin{align}
a_i  + \sum_{i' : i'i\in\mL } s_{i'}\bar{p}_{i'i}   & < s_i  ,&& \text{for each } i\in \mI,\label{RTBP:Stab1} \\
\sum_{\sigma\in\mS_j}  \rho^j_{\sigma}  & \leq 1-\frac{L}{T}, && \text{for } j\in\mJ, \label{RTBP:Stab2} \\
s_i \leq \sum_{\sigma\in\mS_j} & \rho^j_{\sigma}\sigma_i  , && \text{for } j\in\mJ, \label{RTBP:Stab3}
\end{align}}
\hl{where equation \eqref{RTBP:Stab1} represents the need for the
accumulated arrival rates to be less than the potential departure
rates, equation \eqref{RTBP:Stab2} guarantees \hlblue{the yellow and all-red periods at} each cycle maintains
sufficient time for switching and setup between phases\hlblue{,} and equation
\eqref{RTBP:Stab3} indicates the departure rates do not exceed the
allocated service rates. We }let $\bar{\mA}$ denote the closure of
the stability region, that is the set of rates $a=(a_i: i\in\mI)\geq
0$ where the above inequalities in
(\ref{RTBP:Stab1})--(\ref{RTBP:Stab2}) may hold with equality. We
also note that the random variables $A_i(T)$ and the assigned
service time proportions $P^j_{\sigma}(t)$ are corresponding to
$a_i$ and $\rho^j_{\sigma}$ and take their respective values from
the sets $a$ and $\rho$ in the stable case.

Given the vector of queue sizes $(Q_i(t):i\in\mI)$, we define the
total queue size of the road network to be
\begin{equation}\label{QSigma}
Q^\Sigma(t) = \sum_{i\in\mI} Q_i(t).
\end{equation}
So $Q^\Sigma(t)$ gives the total number of cars within the road
network. We say that a policy \hlblue{$P_{\sigma}^{j}(t)$} for serving cars at the junctions
 \hl{stabilizes the network for a vector of
arrival rates $(a_i : i\in\mI)$} if the long run \hlblue{average} number of cars in
the queueing network is finite, in particular,
\begin{equation}\label{stable}
\lim_{T\rightarrow \infty}\bE \left[ \frac{1}{T} \sum_{t=1}^T Q^\Sigma(t) \right] < \infty.
\end{equation}
\hl{This notion of stability originates from the theory of
\hlblue{Markov chains}, where \eqref{stable} gives a necessary and sufficient
condition for positive recurrence, for instance, see
\cite{meyn2009markov}. Our model does not assume \hlblue{that} the underlying
system is \hlblue{Markovian}, thus recurrence cannot be defined. However by
\eqref{stable} we can have the same understanding of necessary and
sufficient conditions  for stability as in the previous literature,
see  \cite{tassiulas1992stability, Br08,sipahi2009stability}. So }in the long run we
expect there to be a finite number of cars within the road traffic
network. If the road network was unstable then we would expect the
number of cars within the system to grow over time. Thus we say that
a policy is unstable \hl{for a vector of arrival rates $(a_i :
i\in\mI)$ if}
\begin{equation}
\lim_{T\rightarrow \infty}\bE \left[ \frac{1}{T} \sum_{t=1}^T Q^\Sigma(t) \right] = \infty.
\end{equation}
We note that if the queue size process was a Markov chain then
definition of stable would be equivalent to the the definition of
positive recurrence for that Markov chain. However, the process that
we will define need not be a Markov chain hence we use the above
definition.

\subsection{Main Theoretical Results}

First of all we state the following known result about the stability
region $\mA$ \hl{defined by
\cref{RTBP:Stab1,RTBP:Stab2,RTBP:Stab3}.} \hlblue{We now refer to
the demand induced by the arrival rates as the load on the network.
In particular, we show that any set of arrival rates outside the
stability region must be unstable no matter what policy is used.
Note that in practice, the traffic load is determined by an
origin-destination (O-D) demand rather than a per-inroad arrival
rate and turning fraction. If the O-D demand is stationary, then
these quantities are also stationary, and the model correctly
captures the load on each roach, and hence the stability of the
network. If the O-D demand is non-stationary, then we capture the
first-order effects by allowing $A_i(t)$ to vary, but the assumption
that $\bar{p}_{ii'}$ is constant is an additional modelling
approximation.}

\begin{proposition}\label{prop1}
Given that the arrivals at each time, $\{A(t)\}_{t=1}^\infty$, are independent
identically distributed random variables with expectation \hl{${a}$, it follows
that
if ${a}\notin \bar{\mA}$ then any policy is unstable under these arrival rates, ${a}$.}
\end{proposition}
The previous proposition shows that the best a policy can do to
stabilize the road traffic network is to be stable for all rates in
$\mA$. The following result shows that our policy is indeed stable
for all arrival rates within the set $\mA$.

\begin{theorem}\label{thrm1}
Given that there exists an $\epsilon>0$ such that for each traffic cycle $t$,
$\bar{a}(t)+\epsilon \mathbf{1} \in \mA$ then, for a constant $K>0$, the long run average queue sizes of in-roads are bounded as\hl{
\begin{equation}\label{Thrm1Eqn}
\lim_{\tau\rightarrow\infty} \bE \left[ \frac{1}{\tau} \sum_{t=0}^{\tau-1}  Q^\Sigma(t)
\right]   \leq   \frac{K}{\epsilon}
\end{equation}}
and thus the policy is stable.
\end{theorem}

We leave the proof of these statements to \ref{sec:apx.proof}. These
results in a certain sense mean, that we can interpret our policy as
being stable for the largest possible set of arrival rates. Thus the
policy provides sufficient throughput in congested traffic as long
as it is possible, reaching an efficient utilization of the existing
capacities. \hlblue{Note that Theorem 1 applies to time varying traffic levels.  Although the stability region $\mathcal{A}$ is multidimensional, the intuition behind the traffic model can be understood by considering the scalar case.  In that case, it corresponds to the expected number of arrivals in each time cycle ($T$) being bounded above.  If we interpret that bound as the traffic level during peak hour, the theorem applies to networks in which queues would remain stable even if peak hour extended indefinitely.  We acknowledge that this is a stricter requirement than necessary, since the system can be stable in the long term even if queues build up during the peak hour, provided they empty sufficiently after the peak.}

\section{Numerical Results - Performance Evaluation and Design}\label{sec:num}

\subsection{Simulation settings}
In this section we evaluate \hl{via simulation} the performance of
our proposed \hl{cyclic phase} BackPressure traffic signal control
and compare its performance with a number of existing self-control
(i.e. decentralized) schemes by \cite{Lammer2010}, \cite{Smith1980}
and \cite{Wongpiromsarn2012} \hl{as detailed below.}

\hl{First, the self-control scheme in \cite{Lammer2010}} aims at
minimizing the waiting times at each intersection anticipating
future arrivals into those intersections instead of just efficiently
clearing exiting queues as in a conventional $\mu c$ priority rule
\cite{Smith1980}.\hl{ However, since future traffic demand is not
known, this scheme~\cite{Lammer2010} more or less \emph{greedily}
attempts to minimize the waiting time. When the setup time or the
amber traffic signal is ignored, this policy (referred to as greedy
policy below)} tends to allocate service to the phase that has
longer queue length. \hl{In contrast, backpressure (includes the one
proposed in this paper) is non-greedy. Backpressure based policy
ensures that an action at this time is not too suboptimal,
regardless of what future traffic is like. Although it does not
explicitly seek to minimize the waiting time, it is likely to result
in lower waiting time and subsequently lower total travel time
through a network than a greedy algorithm that does.}

\hl{Second, the priority rule of the self-control scheme in
\cite{Smith1980} is approximately giving a green time split
proportionally to the total number of vehicles on the in-roads and
thus will be referred to as proportional scheme in the rest of this
section.}

\hl{The third and final} policy in \cite{Wongpiromsarn2012}
allocates green time to the phase that has the highest queue backlog
differences between upstream queues and downstream queue, thus it
will be \hl{referred to as BackPressure policy in this section. Note
that although these benchmarks \hlblue{are in their genesis} by the
standards of currently implemented centralized schemes, they are
state of the art among distributed schemes, and thus appropriate for
the goal of this paper. In summary the following policies will be
evaluated and compared in this section.}

\begin{itemize}
\item\hl{ Cyclic Phase BackPressure policy proposed in this paper: Refer to
Subsection \ref{sec:fixed.BP} for details.}

\item
\hl{BackPressure policy~\cite{Wongpiromsarn2012}:}
\begin{enumerate}
\item
At the beginning of each time slot, based on recent occurrences,
form an estimate of the turning fractions according to
\begin{equation}
\bar{q}_{ii'}(t)= \frac{1}{k} \sum_{\kappa=1}^k p_{ii'}(t-\kappa)
\end{equation}
where $k$ is a parameter of the model.
\item
For each junction $j\in\mJ$, calculate the weight associated with
each service phase at the junction as
 \begin{equation}\label{weightdef}
w_{\sigma}(Q(t)) =  \sum_{i\in j} \sigma_i^j \bigg(  Q_i(t) -
\sum_{i': ii'\in\mL} \bar{q}_{ii'}(t)  Q_{i'}(t)\bigg).
\end{equation}
\item Given these weights, assign the whole service time of the next time slot to phase
$\sigma^*\in\mS_j$ where $w_{\sigma^*}>w_{\sigma}\,\,\forall \sigma
\in \mS_j$.
\end{enumerate}

\item
\hl{Proportional policy~\cite{Smith1980}:}
\begin{enumerate}
\item
At the beginning of each traffic cycle, calculate the weight
associated with each service phase at each junction $j\in\mJ$ as
 \begin{equation}\label{weightdef001}
w_{\sigma}(Q(t)) =  \sum_{i\in \sigma}  Q_i(t) .
\end{equation}
\item Given these weights, within the next service cycle assign an amount of time to each phase $\sigma\in\mS_j$ that is proportional to
\begin{equation*}
P_\sigma^j(t)= \frac{w_{\sigma}(Q(t))}{\sum_{\pi\in \mS_j}
w_{\pi}(Q(t))}.
\end{equation*}
\end{enumerate}

\item
\hl{Greedy policy~\cite{Lammer2010} with no switching time:}
\begin{enumerate}
\item
At the beginning of each time slot, calculate the weight associated
with each service phase at each junction $j\in\mJ$ according to
equation (\ref{weightdef001}).
\item Given these weights, assign the whole service time of the next time slot to phase
$\sigma^*\in\mS_j$ where $w_{\sigma^*}>w_{\sigma}\,\,\forall \sigma
\in \mS_j$.
\end{enumerate}

\end{itemize}

\begin{figure}[h!]
  \centering{
      \includegraphics[width=0.8\textwidth]{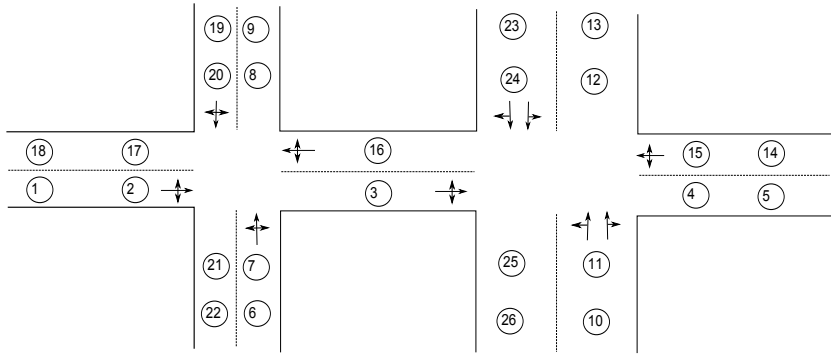}
      \caption{Small network topology.}
  \label{fig:2x1network}
  }
\end{figure}


\begin{figure}[h!]
  \centering{
      \includegraphics[width=1.0\textwidth]{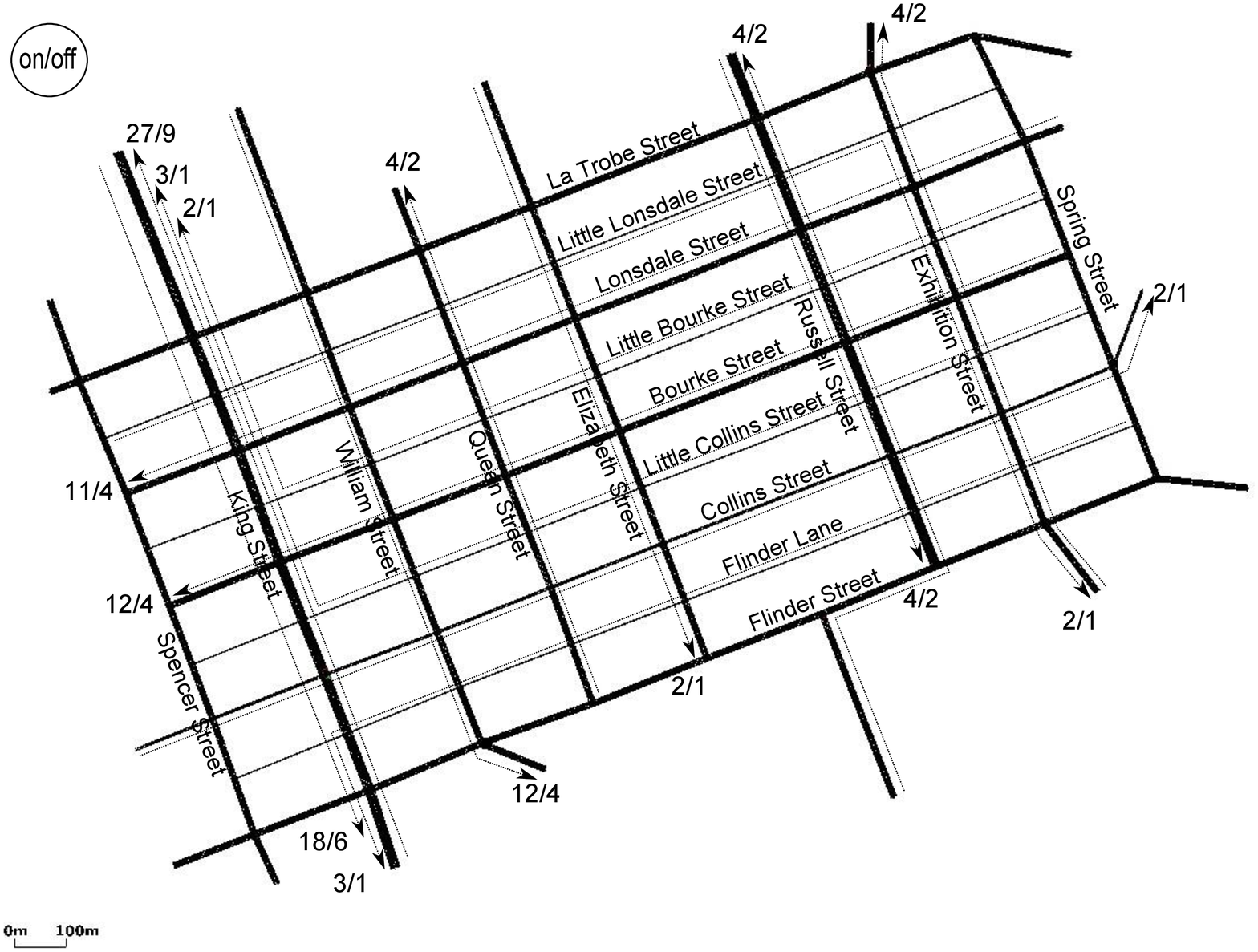}
      \caption{Large CBD network with demands.}
  \label{fig:largeNWdemand}
  }
\end{figure}

\hl{In this section }we utilize an open source microscopic
simulation package SUMO (Simulation of Urban MObility) \cite{SUMO}
to study the above schemes in a small network of two
intersections(Fig.~\ref{fig:2x1network}) and in a large
network(Fig.~\ref{fig:largeNWdemand}) that reassembles the Melbourne
CBD \hl{(Australia)} with about 70 \hl{intersections.}

\hl{The} small network has 2 junctions consisting of several in-roads
(numbered from 1 to 26 on the figure~Fig.~\ref{fig:2x1network}). All
the roads have bi-directional traffic with the North-South road
going through the right intersection having double lanes. Direction
of traffic movements on this network is indicated on each in-road
leading to the junctions. The ingress queues, where vehicles enter
the network, are assumed to be infinite and represented by a set of
long links (i.e. links $\{1, 18, 9, 19, 13, 23, 5, 14, 10, 26, 6,
22\}$ on Fig.~\ref{fig:2x1network}). The cars immediately appear on
the connecting in-roads inside the network, which are of finite
capacity. Since the ingress queues are infinitely large, vehicles
can enter the network even when there is a heavy congestion on the
bottleneck link. All other links (i.e. links$\{2, 3, 4, 7, 8, 11,
12, 20, 21, 24, 25\}$ on Fig.~\ref{fig:2x1network}) have the same
length at $375$ meters which can accommodate maximum $50$ cars per
lane.

The topology of the large CBD Melbourne network is shown in
Fig.~\ref{fig:largeNWdemand}. It consists of $73$ intersections and
$266$ links. Most of the roads are bi-directional except for Little
Lonsdale Street, Little Bourke Street, Little Collins Street and
Flinder Lane which only have a single lane mono-directional traffic.
King Street and Russell Street are the biggest roads in this
scenario, each is modeled as $3$ lanes each direction. Collins
Street has one lane each direction. All other roads have two lanes
each direction. The link lengths are varied between $106$ meters for
the vertical links and $214$ meters and $447$ meters for the
horizontal links except for the ingress links at the edges.

Results are given in terms of the total number of vehicles in the
network and the congestion level which is the average number of
congested links in large network after long simulation runs using
the different control schemes. \hl{In all the studied scenarios, the
exact queue lengths and turning fractions are observed directly from
the simulation and used to make control decision in various
policies.} These variables are calculated using Matlab \cite{MATLAB}
based on the actual control algorithm and then are fed back into the
SUMO simulation at every time step. We ignore switching times (i.e.
transition between phases) in all control schemes in our study. This
overhead can be incorporated into the simulation by extending the
phase times. Nevertheless, the qualitative insights gained in this
section would not change by that extension.

\subsection{Performance Study}
\hl{Below} we evaluate the performance of our \hl{cyclic phase
BackPressure} scheme and compare it with other policies using fixed
setting of routes in the studied networks \hl{using simulation}.
\hl{The cycle time of the cyclic phase BackPressure policy and the
proportional policy were set to $30$ seconds, while the slot time of
the BackPressure policy and the greedy policy were set to $10$
seconds in our simulation.}

\subsubsection{Small Network}
\begin{figure}[h!]
  \centering{
      \includegraphics[width=0.8\textwidth]{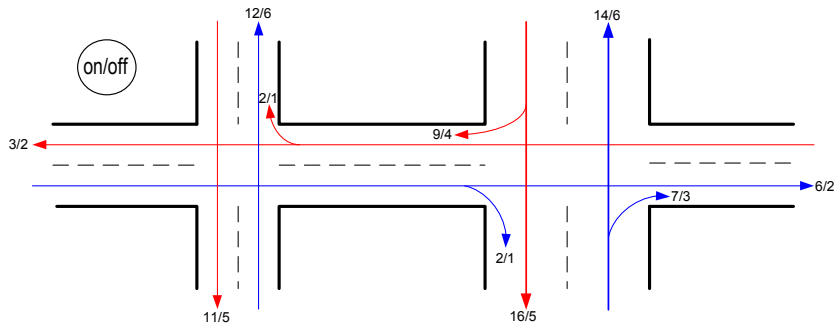}
      \caption{Small network with demands.}
  \label{fig:case2}
  }
\end{figure}

First, we study a small network scenario, for which the turning
information and the arrival rates are indicated in
Fig.~\ref{fig:case2}. In particular, the arrows indicate 10 routes
with direction and demands (cars/minutes) in the peak and off-peak
(i.e. on/off) time periods as shown in Fig.~\ref{fig:case2}. The
$\eta$ value introduced in (\ref{policydef}) was set to $2.5$. The
main traffic flows are the ones with North-South direction of the
second junction. The two BackPressure policies give the majority of
service time to the North-South phase of that junction which leads
to heavy congestion on links 1, 2, 3, 14 and 15. On the other hand,
the proportional policy and greedy policy put more balance between
the service times depending on the queue lengths which creates more
congestion in the North-South direction at the cost of having less
congestion in East-West direction.

\begin{figure}[h!]
  \centering{
      \includegraphics[width=0.8\textwidth]{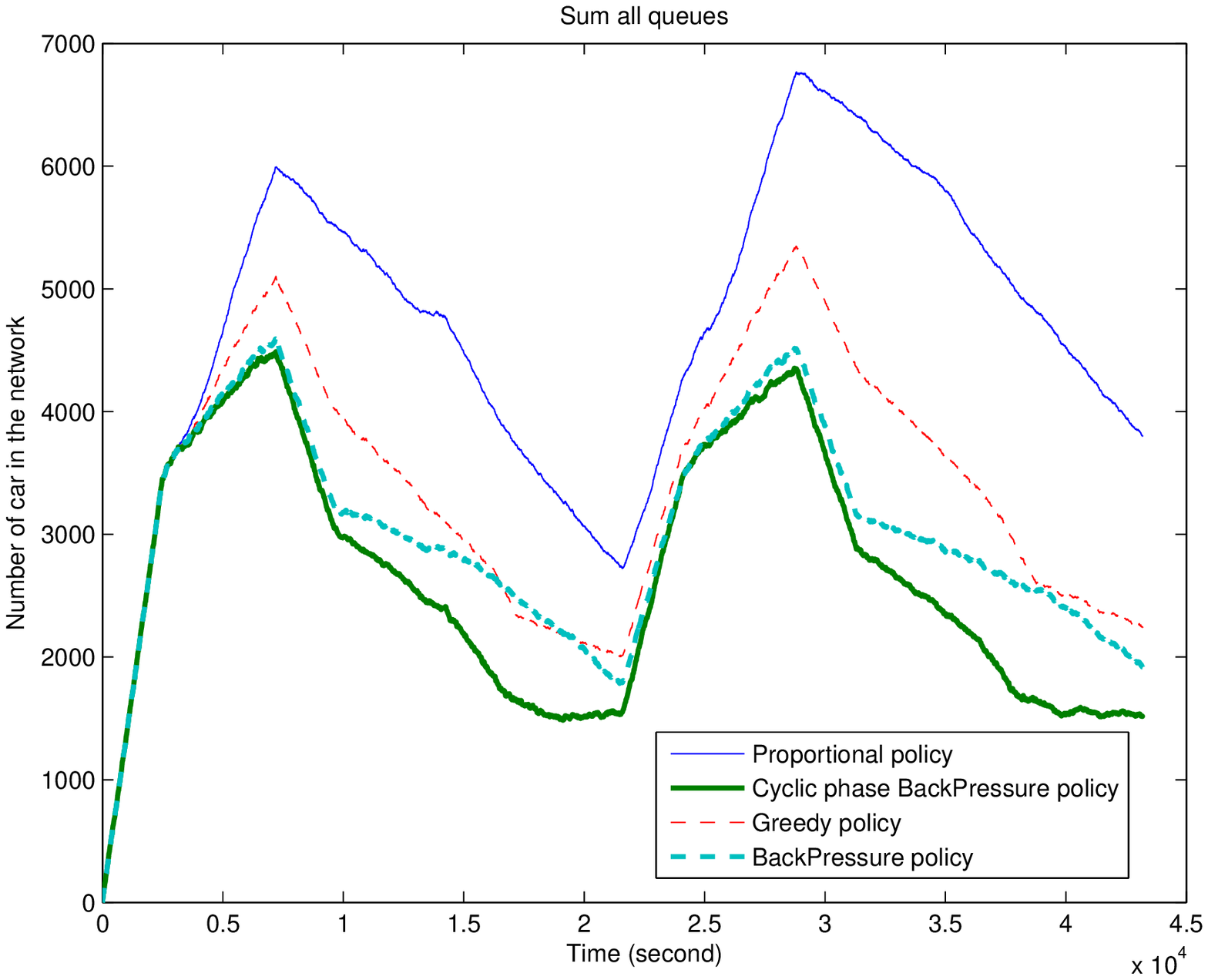}
      \caption{Throughput for the small network using different policies. Cycle time for the cyclic phase BackPressure policy and the proportional policy is 30 seconds. Slot time for the BackPressure policy and the greedy policy is 10 seconds.}
  \label{fig:smallNWsumall}
  }
\end{figure}

\begin{figure}[h!]
  \centering{
     \includegraphics[width=0.8\textwidth]{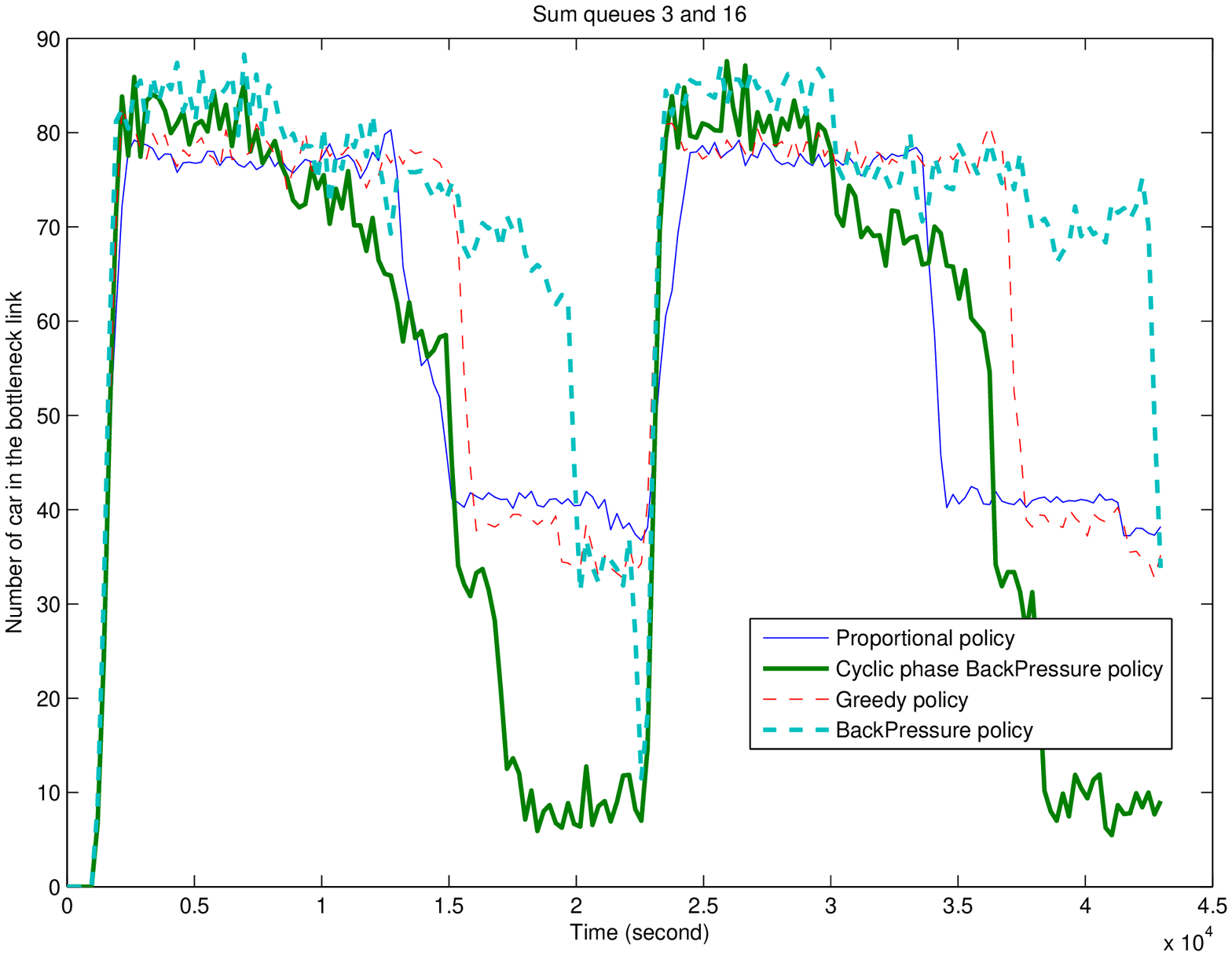}
      \caption{Congestion of bottleneck link for the small network using different policies. Cycle time for the cyclic phase BackPressure policy and the proportional policy is 30 seconds. Slot time for the BackPressure policy and the greedy policy is 10 seconds.}
  \label{fig:smallNWcongestion}
  }
\end{figure}

Results are shown in Fig.~\ref{fig:smallNWsumall} and
Fig.~\ref{fig:smallNWcongestion} where the total number of vehicles
in the network and in the congested link between the two junctions
are plotted over time. Note that there are two in-roads between the
two intersections but only link 3 is congested due to large traffic
flows in the North-South direction at the second intersection.

\hl{Note that Figs \ref{fig:smallNWsumall} and
\ref{fig:smallNWcongestion} were based on the number of cars present
at the times when control decisions were made which is $30$ seconds
for cyclic phase BackPressure and proportional policies. In
contrast, the average travel time depends on the waiting time on
individual link which is an integral of queue size over continuous
time. For this reason, intermediate queue size was also measured at
10\,s intervals in the simulation, and the results differed by less
than \hlblue{2\% in compare with the coarse sampling at once per cycle assuming linear interpolation}. The resulting travel time values are reported in the next
sub-section~\ref{sub-sec:Design}.}

Observe that the \hl{cyclic phase BackPressure} control yields a lower
number of total vehicles present in the network and thus results in
higher number of vehicles reaching their destination (i.e. increased
network throughput) during the whole simulation. This is due to the
fact that in the \hl{two BackPressure control schemes} when the
bottleneck link (link 3) is congested, less green time will be
allocated to the East-West direction at the first junction. As a
result more traffic can move through the North-South direction and
the impact of a spill back from the second junction on the overall
network throughput decreases. Similarly, the \hl{BackPressure policy}
also outperforms the \hl{proportional policy and the greedy policy,}
since those two schemes allocate similar amount of green time to the
East-West direction at the first junction despite the presence of a
spilled back traffic and thus waste some of the green time.

\subsubsection{Large Melbourne CBD network}
A similar study is performed with a large network with its turning
information and arrival rates indicated in
Fig.~\ref{fig:largeNWdemand}. The parameter $\eta$ is once again set
to $2.5$. In this setting, the King Street has the largest flows,
thus, any flow that shares an intersection with King Street tends to
be under-served especially the intersection between King Street and
Lonsdale Street and the intersection between King Street and Bourke
Street. Generally in the peak period, congestion in any link will
cause spill-back which leads to further congestion in the
neighbouring junctions. This can only be recognized by the \hl{two
BackPressure} policies through comparing the in-road $i$ and
out-road $i'$, and more service will be allocated in this case to
traffic flows on the less congested directions. In contrast, the
\hl{proportional policy and the greedy policy} only consider the
queue lengths present at in-road $i$, and may waste some green time
to the congested direction where traffic comes to a standstill due
to the spill back.

\begin{figure}[h!]
  \centering{
      \includegraphics[width=0.8\textwidth]{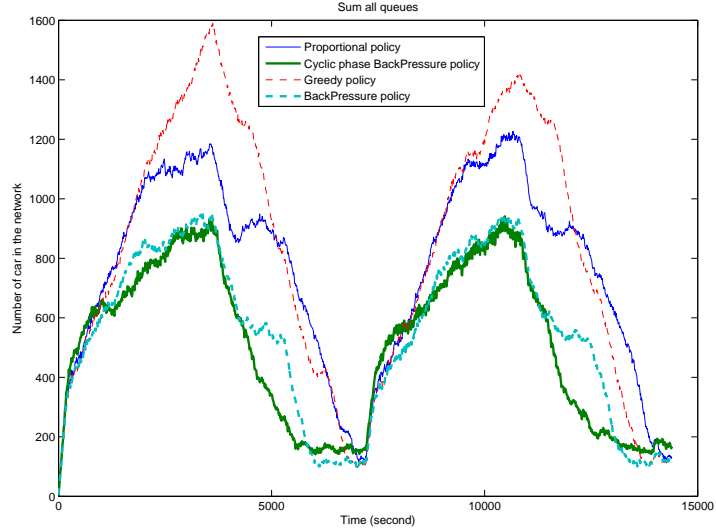}
      \caption{Throughput for the large network using different policies. Cycle time for the cyclic phase BackPressure policy and the proportional policy is 30 seconds. Slot time for the BackPressure policy and the greedy policy is 10 seconds.}
  \label{fig:largeNWsumall}
  }
\end{figure}
\begin{figure}[h!]
  \centering{
      \includegraphics[width=\textwidth]{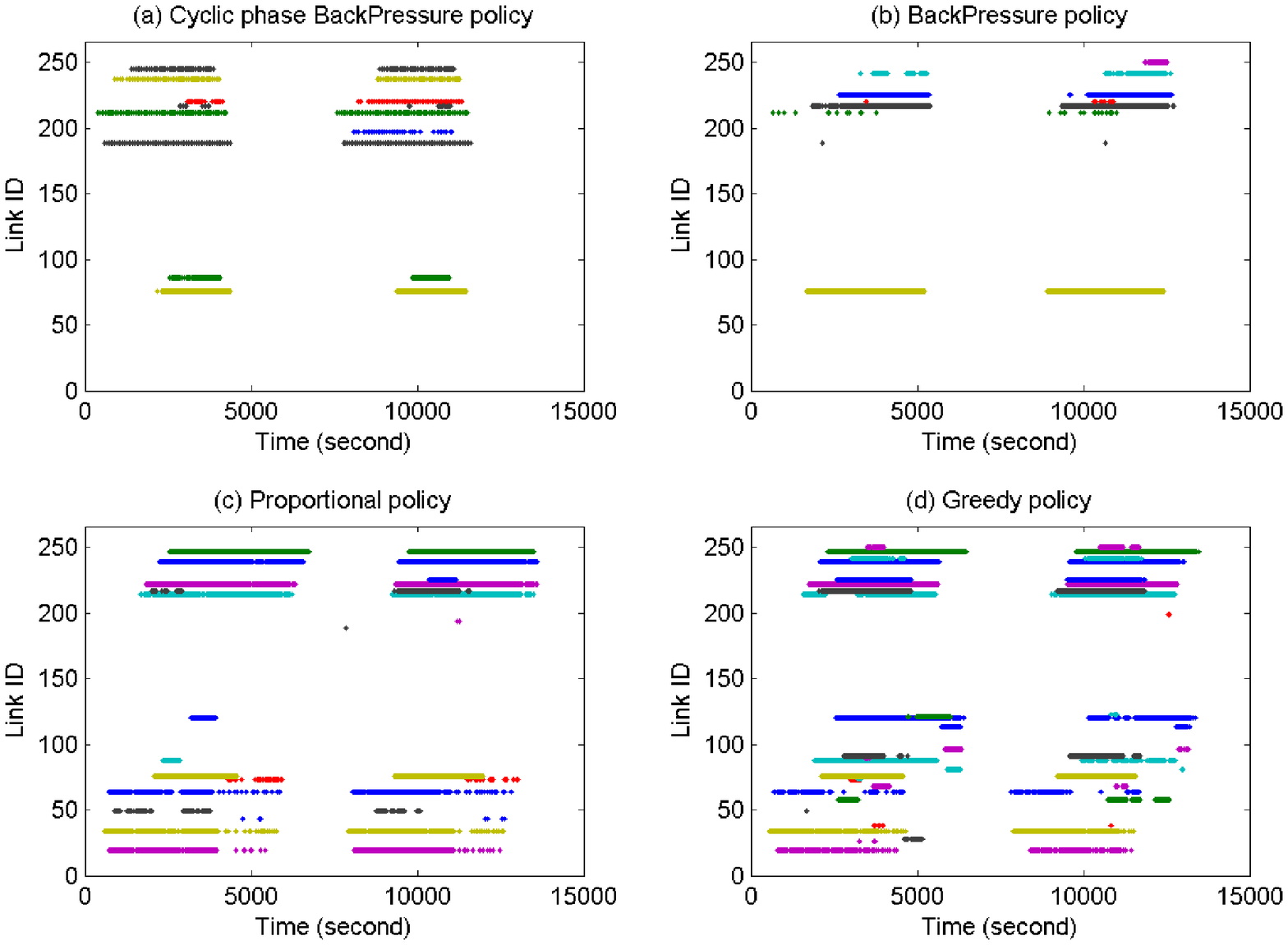}
\caption{Congestion level for the large network using different
policies. Cycle time for the cyclic phase BackPressure policy and
the proportional policy is 30 seconds. Slot time for the
BackPressure policy and the greedy policy is 10 seconds.}
  \label{fig:largeNWcongestion}
  }
\end{figure}

Results for this scenario are shown in Fig.~\ref{fig:largeNWsumall}
and Fig.~\ref{fig:largeNWcongestion}. As shown in
Fig.~\ref{fig:largeNWsumall} \hl{the cyclic phase BackPressure policy} has
the lowest total number of vehicles in the network, thus it provides
the highest throughput, whereas the second highest throughput is
provided by the \hl{BackPressure policy}. Similarly to
the previous scenario, the two BackPressure policies outperform the
other two policies in case of heavy congestion because they take
into account downstream queue lengths, and thus allocate resources
(i.e. service phases) more efficiently.

Fig.~\ref{fig:largeNWcongestion} plots the congested link over time.
\hl{Herein a link} is said to be congested at a certain time if its queue
length is more than $85\%$ of the link capacity. It is clear that
the \hl{two BackPressure policies} reduce the number of congested links
significantly (i.e. less number of vehicles inside the network)
resulting in higher network throughput.

\subsection{Experimental Parameter Design} \label{sub-sec:Design}
The cycle length in the \hl{cyclic phase BackPressure policy and the
proportional policy} and the frequency of making decision in the
\hl{BackPressure policy and the greedy policy} play a crucial role
in the performance of the control scheme. A long cycle length or the
low frequency of making decision may be less efficient due to the
fact that the queue might be depleted before the end of the service
time. In the other hand, a short cycle length may reduce the overall
capacity since the vehicles have to stop and accelerate more often.
Note that the latter is in fact represents a switching cost between
phases even though the amber traffic signal is not considered here.
This subsection investigates the impact of the cycle time and
decision making frequency on the throughput and congestion level of
each scheme. We study both network topologies (the small network and
the large CBD network) under the similar demand levels as in the
previous subsection with different cycle times and decision
frequencies. Particularly, for the \hl{cyclic phase BackPressure
policy and the proportional policy}, the cycle length is set to
$\{30,60,90,120\}$ seconds, and for the \hl{BackPressure policy and
the greedy policy}, a decision is made every $\{10,30,60,90\}$
seconds, respectively.

\subsubsection{Small Network}
\begin{figure}[h!]
  \centering{
      \includegraphics[width=0.8\textwidth]{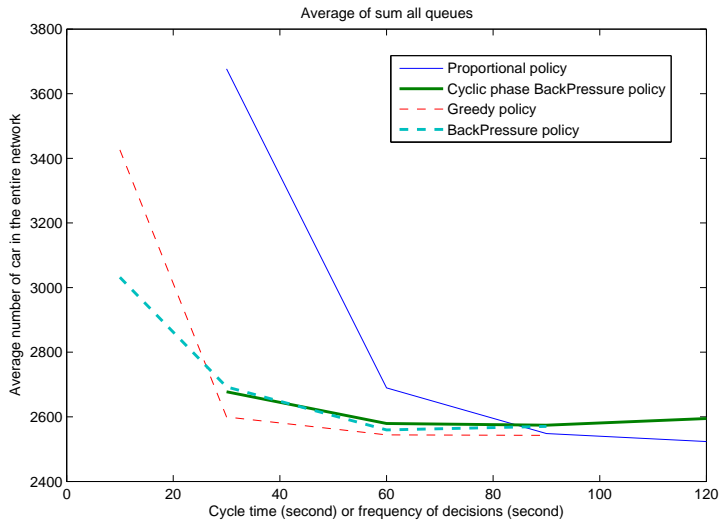}
      \caption{Throughput vs cycle time/frequency of decision for the small network.}
  \label{fig:smallNWdesignthroughput}
  }
\end{figure}

\begin{figure}[h!]
  \centering{
      \includegraphics[width=\textwidth]{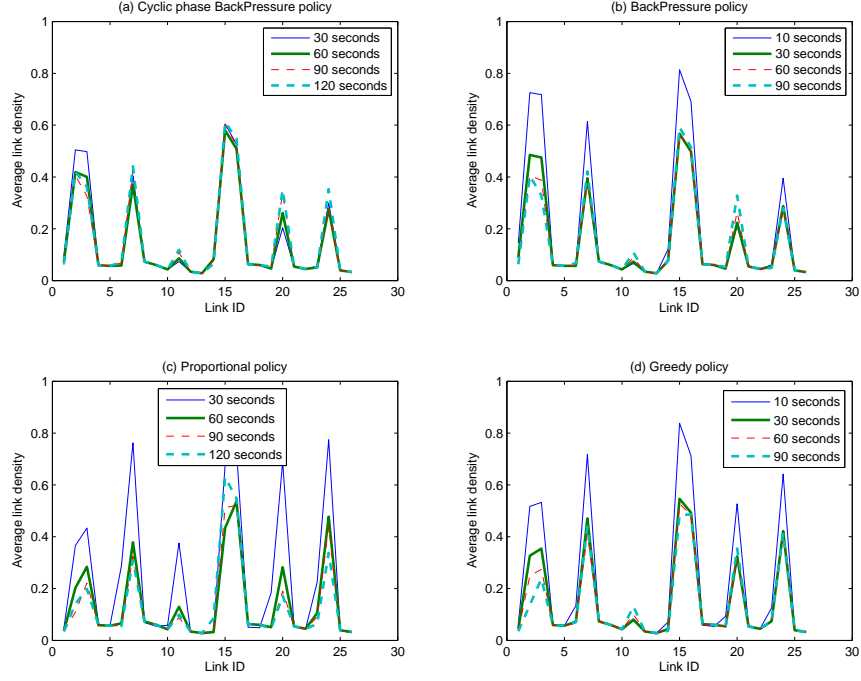}
      \caption{Average link density for the small network.}
  \label{fig:smallNW_avg_density}
  }
\end{figure}

\begin{figure}[h!]
\centering{
      \includegraphics[width=\textwidth]{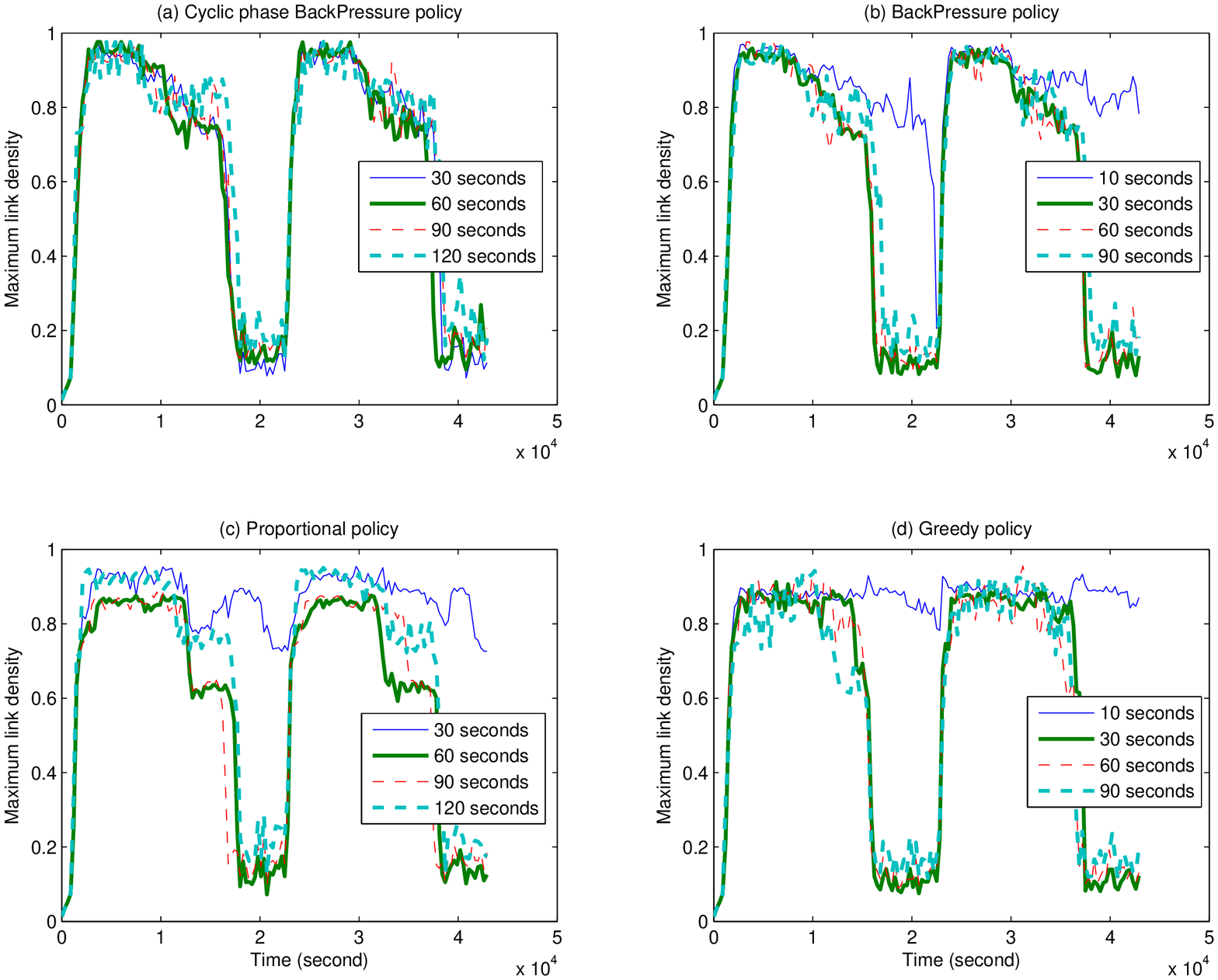}
      \caption{Maximum link density for the small network.}
  \label{fig:smallNW_max_density}
}
\end{figure}

For small network, the results are presented in
Fig.~\ref{fig:smallNWdesignthroughput},
Fig.~\ref{fig:smallNW_avg_density}, and
Fig.~\ref{fig:smallNW_max_density}.
Fig.~\ref{fig:smallNWdesignthroughput} shows the average number of
vehicles in the network plotted against different cycle times.
Because vehicle does not disappear and stays in the network until it
exists, lower number of vehicles in the network equates to higher
throughput of the same demand. In this scenario, all of the studied
policies provide a similar throughput using their corresponding best
setting. Furthermore, it can be observed that in congested network
higher cycle times tend to have better throughput because the
traffic flows are less interrupted by the switching between phases.
Nevertheless, the proposed cyclic phase BackPressure is less
sensitive to the changes of cycle time while producing a compatible
throughput.

Fig.~\ref{fig:smallNW_avg_density} plots the average link densities
versus link ID. It shows that congestions occur in the same set of
links throughout all of the studied policies. In overall, \hl{the
cyclic phase BackPressure policy and the proportional policy} have
lower link density than the other two policies.
Fig.~\ref{fig:smallNW_max_density} shows the maximum link density
over time pointing to when the congestions occur in the network. In
all cases, congestions appear during the peak period and some
portion of time during the off-peak period before the build-up
traffic can be substantially drained. However, when the cycle time
or the decision frequency is set too small (e.g. 10 seconds for the
\hl{BackPressure policy and the greedy policy} and 30 seconds for
the \hl{cyclic phase BackPressure policy and the proportional
policy}), the congestions are not able to cleared at all except for
our cyclic phase BackPressure policy.


\subsubsection{Large Melbourne CBD Network}
\begin{figure}[h!]
  \centering{
      \includegraphics[width=0.8\textwidth]{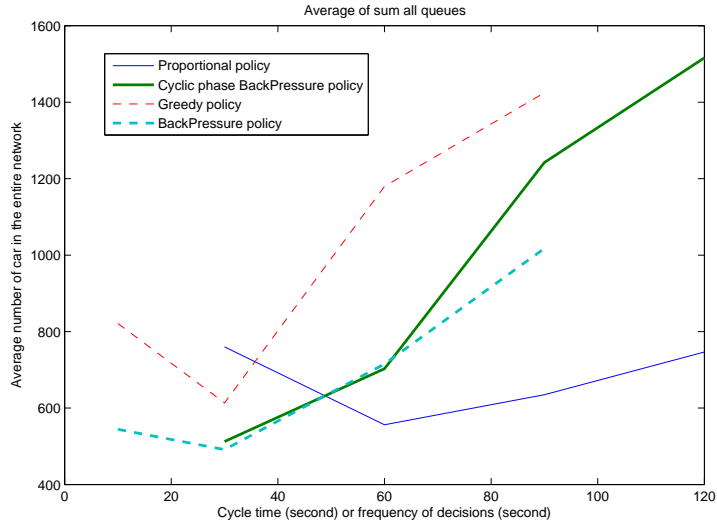}
      \caption{Throughput vs cycle time/frequency of decision for the large network.}
  \label{fig:largeNWdesignthroughput}
  }
\end{figure}

\begin{figure}[h!]
\centering{
      \includegraphics[width=\textwidth]{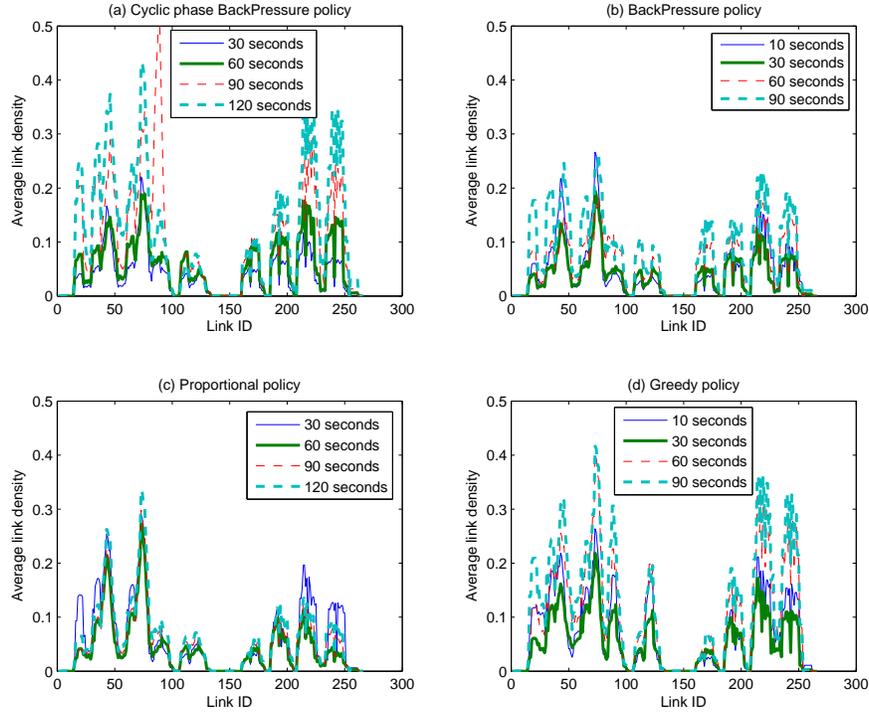}
\caption{Average link density for the large network.}
\label{fig:largeNW_avg_density} }
\end{figure}

\begin{figure}[h!]
\centering{
      \includegraphics[width=\textwidth]{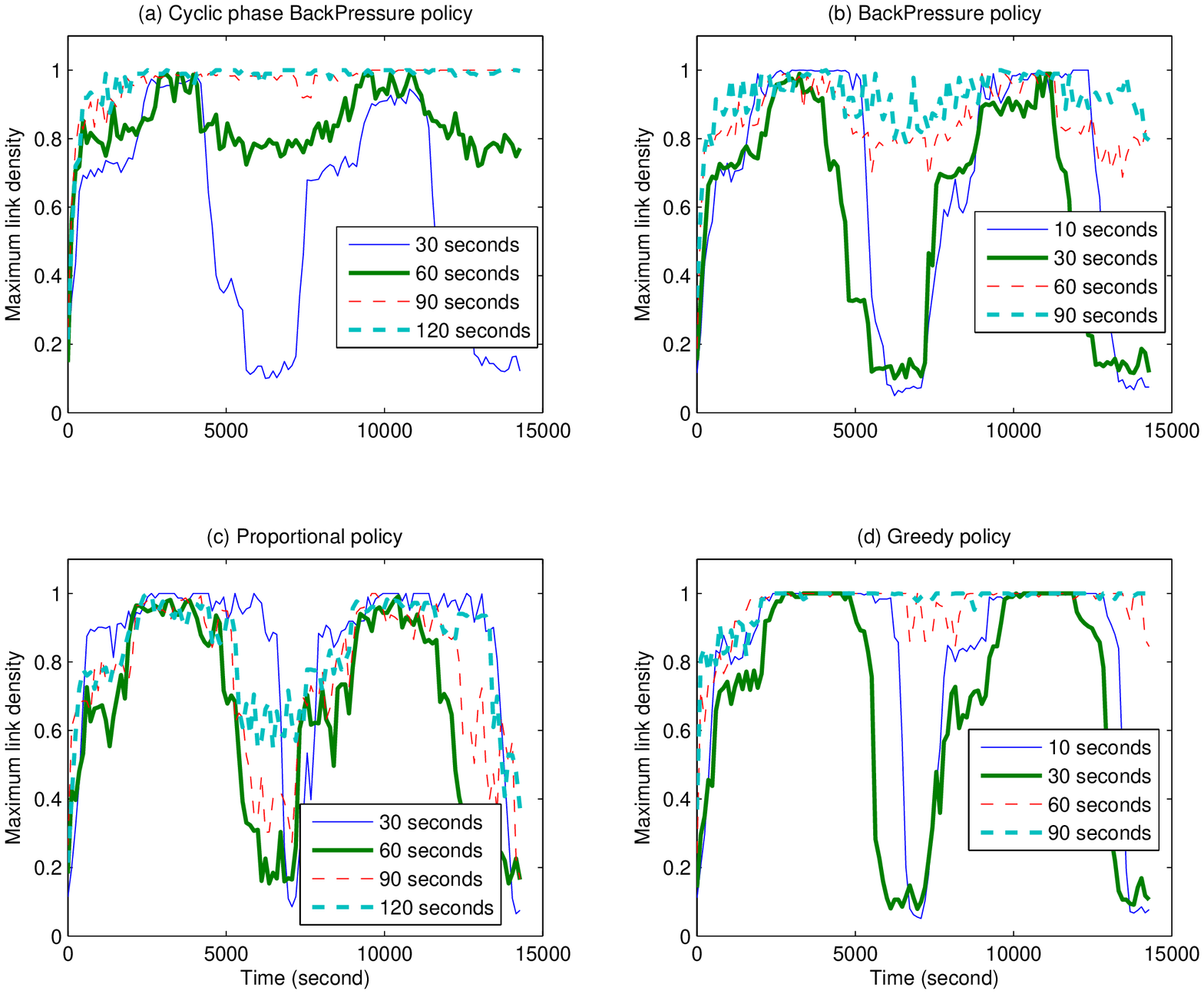}
\caption{Maximum link density for the large network.}
\label{fig:largeNW_max_density} }
\end{figure}

The impact of the cycle time and frequency of decision making on
network throughput and congestion level using different policies for
a large network are investigated and discussed in this subsection.
The results are shown in Fig.~\ref{fig:largeNWdesignthroughput},
Fig.~\ref{fig:largeNW_avg_density} and
Fig.~\ref{fig:largeNW_max_density}.

The average number of vehicles in the network for each setting is
presented in Fig.~\ref{fig:largeNWdesignthroughput}. Unlike the
results in the small network, there clearly exists an optimal value
for the cycle length or decision making frequency of each policy. In
particular the optimal cycle length for \hl{proportional
policy} is $60$ seconds, while the optimal cycle length/decision
frequency for all other policies is $30$ seconds.


\begin{table}[h!]
\centering{
\begin{tabular}{ c | c  c c c }
  & Proportional & Cyclic phase BP & Greedy & BP \\ \hline \hline
  Avg. travel time & 478.0 & 409.6 & 514.0 & 408.5\\
\end{tabular}
\caption{Average travel time (in seconds) for the Melbourne CBD
network using optimal setting for each policy (i.e. 60 seconds cycle
length for the proportional policy, and 30 seconds cycle length for
all other policies).} \label{table:avgTravelTime} }
\end{table}

Observe that the optimal cycle length in the large network scenario
is shorter than that of the small network, which can be explained by
shorter link lengths and larger number of intersections on any
route. Both increase the interdependency between intersections and
their performance as arriving traffic into any internal intersection
is an output traffic from the others.

\hl{Furthermore, the average travel time of each policy in their
best setting of cycle length and decision making frequency is shown
in Table \ref{table:avgTravelTime}. There is a strong correlation
between the average number of vehicles in the network and average
travel time through that network. In particular, the higher number
of vehicles in the network results in the longer travel time and
vice versa. It can be seen that the proposed cyclic phase
BackPressure policy has a competitive average travel time between
all the considered policies, and the results show that the
BackPressure-based policies yield a significant\hlblue{ly} better average
travel time \hlblue{than} that of the greedy policy. Note that this
better average travel time has been achieved with the control
decisions using the queue size measurements at discrete time
intervals (once in every cycle) only as explained earlier.}

Fig.~\ref{fig:largeNW_avg_density} plots the average link densities
against the link ID. It shows that the congestion area is varied
with different policy and with different parameter settings. Any
cycle length/decision frequency setting other than the optimal
setting obviously increases the congestion greatly.

Finally, Fig.~\ref{fig:largeNW_max_density} shows the maximum link
density over time. The cycle length plays a vital role to prevent
congestion in this scenario. In the \hl{cyclic phase BackPressure policy,}
the $30$ second cycle length is undoubtedly outstanding. In other
policies, the small cycle lengths are seemed to be better due to the
short link lengths.
%

\section{Conclusion}
\label{sec:concl}

We proposed in this paper a novel decentralized signal control
strategy based on the so-called BackPressure policy that does not
require any $a ~priori$ knowledge of the traffic demand and only
needs information (i.e. queue size) that is local to the
intersection. In contrast to other existing BackPressure-based
policies in which phases can form an erratic and unpredictable order
resulting in potential unsafe operation, our \hl{scheme allocates}
non-zero amount of time to each phases within the cycle, thus
\hl{repeating them in a cyclic manner}. Furthermore, unlike all the
other existing Backpressure-based policies, no knowledge of the
local turn ratios (turning fractions) is required in our control
strategy. Instead any unbiased estimator of the turning fractions
can be utilized in the proposed scheme. We have formally proved the
stability results of the proposed signal control policy even though
the controllers are reacting based only on local information and
demand in an distributed manner. The stability results indicate that
our policy is stable for the largest possible set of arrival rates
(or demand) that will provide sufficient throughput even in
congested network.

Using simulation, we compared our \hl{cyclic phase BackPressure}
performance against other well-known policies in terms of network
throughput and congestion level using both small and large network
topology with fixed routings. The results showed that our \hl{cyclic
phase BackPressure policy} tends to outperform other \hlblue{distributed} polices
both in terms of throughput and congestion. Although the performance
of each policy varies widely depending on the parameter setting such
as cycle length or decision frequency, under the optimal setting
among the cases studied, the \hl{BackPressure with cyclic and
non-cyclic operation} have better throughput in compare with the
other policies.

There are still many issues, such as non-constant switching times,
finite link travel time and link capacity etc. that have not been
considered here and will be a subject of future work.

\section*{Acknowledgements}
This work was supported by the Australian Research Council (ARC)
Future Fellowships grants FT120100723 and FT0991594.

\appendix
\section{Estimation of Turning Fractions}\label{sec:apx.est}

An important aspect not addressed in the previous studies applying
back pressure is the estimation of \hl{traffic turning fractions.}
Previous studies have either assumed the turning fractions are
either explicitly known or have been calculated prior to the
implementation of the policy.

Here we emphasize that the turning fractions can be estimated using
recent locally calculated information about traffic flows. For
instance, if we form an estimate on the turning fractions based on
the last $k$ service cycles
\begin{equation}
\bar{q}_{ii'}(t)= \frac{1}{k} \sum_{\kappa=1}^k \hat{p}_{ii'}(t-{\kappa}),
\end{equation}
\hl{where $ \hat{p}_{ii'}(t)$ denotes the measurement result for
$p_{ii'}(t)$.} If $Q_i(t)=0$ then any estimate may be used to define
$p_{ii\rq{}}(t)$. Given that \hl{turning fractions} are stationary
and independent of queue sizes, these estimates form an unbiased
estimate of the underlying \hl{turning} probabilities \hl{as long as
the measurement error in $\hat{p}_{ii'}(t)$ has zero mean,} since
then
\begin{equation}
\bar{p}_{ii'} = \bE [ \bar{q}_{ii'}(t) | Q_i(t) >0 ].
\end{equation}
Other rules incorporating historical data or more recent data could
also be considered here. What is necessary is that
$\bar{q}_{ii'}(t)$ provides an unbiased estimate of the underlying
turning fractions of the vehicles for non-empty queues. It is even
possible to use an inconsistent estimate of \hl{the turning
fractions or for the proportions to change on a larger time scale,}
as long as the estimate is unbiased, independent of (the history of)
$Q(t)$.

\section{Proof of the main stability result}\label{sec:apx.proof}
In this section we prove Theorem \ref{thrm1}. In order to that, we
have to clarify some assumptions made about the stochastic elements
of our model and provide some technical lemmas which are proven in a
supplementary document together with Proposition \ref{prop1}.

\subsection{Assumptions}\label{sec:apx.assumptions}
\begin{enumerate}
\item The number of cars that can be served from any in-road within a traffic cycle is bounded,
\begin{equation}
S_{\max} = \max_{t\in\bZ_+, i\in\mI} S_i(t) < \infty.
\end{equation}

\item $(p_{ii'}(t): ii'\in\mL)$ is stationary and
independent of queue lengths $(Q_i(\tau):i\in\mI)$ and the number of cars served at each queue
$(S_i(\tau):i\in\mI)$ for all $\tau \leq t$.

\item The matrix $I - \bar{p}$ is invertible. Thus we have, for $ii'\in\mL$
\begin{align}\label{Spmean}
\bE \left[ S_i(t)p_{ii'}(t) | Q(t) \right] &=\sum_{\sigma \in \mS_j} \sigma_i P^j_\sigma(t)\bar{p}_{ii'}\\
\bE \left[ S_i(t)\wedge Q_i(t) p_{ii'}(t) | Q(t) \right] &=\bE \left[ S_i(t)\wedge Q_i(t) | Q(t) \right]\bar{p}_{ii'}.
\end{align}

\item The number of arrival $(A_i(t): t\in \bZ_+^I)$ is independent of the state of the queues in the road traffic network. Thus the average arrival rate into each junction can be defined as
\begin{equation}\label{abar}
\bar{a}_i(t)=\bE[A_i(t)].
\end{equation}

\item \hl{The error term in the queue size measurement, $\delta(t)$ is bounded, i.e.
\begin{equation}\label{error_bound}
|\delta_i(t)|\leq\delta_{\text{max}}\quad\forall i\in\mI.
\end{equation}}
\end{enumerate}

\subsection{Lemmas}\label{sec:apx.lemmas}

 In this section, we will prove a number of additional lemmas that are required for the main proofs.
\hl{The first lemma describes the difference of the weights caused by the error in measurement, whereas the second describes a general result on weights. The third lemma is a consquence of them, introducing a bound which is used multiple times in proving later statements.}
\hl{\begin{lemma}\label{weightdiflemma}
\begin{equation}\label{weightdifbymeas}
 w_{\sigma }(Q(t))= w_{\sigma }(\hat{Q}(t))- w_{\sigma }(\delta(t)),
\end{equation}
where $w_{\sigma }(Q(t))$ is defined according to \eqref{wdef}.
\end{lemma}}
\begin{proof}
This statement is a straightforward consequence of the definition of the weights and the measurement error given in \eqref{hat_error} and \eqref{wdef}. Namely,
\begin{equation}
\begin{aligned}
 w_{\sigma }(Q(t))&=\sum_{i\in j} \sigma_i \left( Q_i(t) - \sum_{i': ii'\in\mL} \bar{q}_{ii'}(t)  Q_{i'}(t)\right)\\
 &=\sum_{i\in j} \sigma_i \left( \left(\hat{Q}_i(t)-\epsilon_i(t)\right) - \sum_{i': ii'\in\mL} \bar{q}_{ii'}(t)  \left(\hat{Q}_{i'}(t)-\delta_{i'}(t)\right)\right)\\
 &=\sum_{i\in j} \sigma_i \left( \hat{Q}_i(t) - \sum_{i': ii'\in\mL} \bar{q}_{ii'}(t)  \hat{Q}_{i'}(t)\right)\\
 &-\sum_{i\in j} \sigma_i \left( \delta_i(t) - \sum_{i': ii'\in\mL} \bar{q}_{ii'}(t)  \delta_{i'}(t)\right)\\
 &= w_{\sigma }(\hat{Q}(t))- w_{\sigma }(\delta(t)).
\end{aligned}
\end{equation}
\end{proof}

\begin{lemma}\label{weightslemma}
Given weights $(w_y :y\in\mY)$ with elements indexed by finite set $\mY$, we consider $Y$ a random variable with the following probability of event $y$:
\begin{equation}
P_y= \frac{e^{\eta w_y}}{\sum_{y'\in\mY}e^{\eta w_{y'}}}
\end{equation}
then, the expected value of the weights $w$ under this distribution obey the following inequality
\begin{equation}
\bE w_Y \geq \max_{y\in\mY} w_y - \frac{1}{\eta}\log |\mY|.
\end{equation}
\end{lemma}
\begin{proof}
In the following inequality,  we note that the entropy of a distribution $H(P)=-\bE \log P(Y)$ is maximized by a uniform distribution on $\mY$, $H(U)=\log |\mY|$.
\begin{align*}
\bE w_Y = &\frac{1}{\eta} \log\left( \sum_{y\in\mY} e^{\eta w_y} \right) + \frac{1}{\eta}\bE \log P(Y) \geq \frac{1}{\eta} \log\left( \sum_{y\in\mY} e^{\eta w_y} \right) -\frac{1}{\eta}\log{|\mY|}\\
\geq &\frac{1}{\eta}   \log \Big( e^{\eta \max_{y\in\mY} w_y}\Big) - \frac{1}{\eta}\log{|\mY|} = \max_{y\in\mY} w_y - \frac{1}{\eta}\log |\mY|,
\end{align*}
as required.
\end{proof}

\begin{lemma} \label{weightslemma2}
\begin{equation}
 \sum_{\sigma \in \mS_{j}}    P^{j}_\sigma  w_{\sigma }\left(\hat{Q}(t)\right)\geq\max_{\sigma\in\mS_j}\Big\{   w_{\sigma}\left(Q(t)\right) \Big\}-\delta_{\max}S_{\max}|\mI|- \frac{1}{\eta} \log |\mS_j|.
\end{equation}
\end{lemma}
\begin{proof}
We can prove this statement by applying Lemma \ref{weightslemma} to $w_{\sigma }(\hat{Q}(t))$ , expanding the terms according to Lemma \ref{weightdiflemma} and using the bound on $\delta(t)$ in \eqref{error_bound}.
\begin{equation}
\begin{aligned}
 \sum_{\sigma \in \mS_{j}}    P^{j}_\sigma  w_{\sigma }(\hat{Q}(t)) &=\bE w_{\sigma }(\hat{Q}(t)) \geq \max_{\sigma\in\mS_j}\Big\{   w_{\sigma}\left(\hat{Q}(t)\right) \Big\}- \frac{1}{\eta} \log |\mS_j|\\
 &\geq\max_{\sigma\in\mS_j}\Big\{   w_{\sigma}\left(Q(t)\right)+w_{\sigma}\left(\delta(t)\right) \Big\}- \frac{1}{\eta} \log |\mS_j|\\
 &\geq\max_{\sigma\in\mS_j}\Big\{   w_{\sigma}\left(Q(t)\right) \Big\}+\min_{\sigma\in\mS_j}\Big\{w_{\sigma}\left(\delta(t)\right)\Big\}- \frac{1}{\eta} \log |\mS_j|\\
  &\geq\max_{\sigma\in\mS_j}\Big\{   w_{\sigma}\left(Q(t)\right) \Big\}-\delta_{\max}S_{\max}|\mI|- \frac{1}{\eta} \log |\mS_j|
\end{aligned}
\end{equation}
\end{proof}
Lemma \ref{QKLemma} and Lemma \ref{QqqLemma} introduce bounds on the increments of the square of the queue sizes and their conditional expectations.
\begin{lemma}\label{QKLemma}
There exists a constant $K_0\geq 0$ such that our queue size process, \eqref{Qprocess}, obeys the  bound
\begin{equation}
\frac{1}{2} Q_i(t+1)^2 - \frac{1}{2} Q_i(t)^2  \leq Q_i(t) \left( A_i(t) - S_i(t)
 +\sum_{i':i'i\in \mL} S_{i'}(t)p_{i'i}(t) \right)   + K_0.
\end{equation}
\end{lemma}
\begin{proof}
Firstly, the following bound holds for the queue size process, \eqref{Qprocess}.
\begin{align}
Q_i(t+1) &=Q_i(t) - S_i(t)\wedge Q_i(t)   + A_i(t) + \sum_{i' : i'i\in\mL }[S_{i'}(t)\wedge Q_{i\rq{}}(t)]   p_{i'i}(t) \notag \\
& \leq Q_i(t) - S_i(t)\wedge Q_i(t)   + A_i(t) + \sum_{i' : i'i\in\mL }S_{i'}(t)  p_{i'i}(t) \notag \\
& \leq
\begin{cases}
Q_i(t) - S_i(t)   + A_i(t) + \sum_{i' : i'i\in\mL }S_{i'}(t)  p_{i'i}(t), &\text{if } Q_i(t) \geq S_i(t),\\
 A_i(t) + \sum_{i' : i'i\in\mL }S_{i'}(t)  p_{i'i}(t), &\text{otherwise}.
\end{cases}\label{Qcases}
\end{align}
Let\rq{}s consider the two cases above. Firstly, if $Q_i(t) < S_i(t)$ then, and according to the above bound, we have
\begin{align}
 \frac{1}{2} Q_i(t+1)^2- \frac{1}{2} Q_i(t)^2 \leq &\frac{1}{2} Q_i(t+1)^2 \leq \frac{1}{2}\left( A_i(t) + \sum_{i' : i'i\in\mL }S_{i'}(t)  p_{i'i}(t) \right)^2\notag  \\
\leq &\frac{1}{2}\left(a_{max} + S_{\max} (|\mI|+1)\right)^2 \notag\\
\leq &- Q_i(t) \left(  S_i(t) - A_i(t) - \sum_{i' : i'i\in\mL }S_{i'}(t)  p_{i'i}(t) \right) \label{Qterm1} \\
& +  S_{\max} \left(a_{max} + S_{\max} (|\mI|+1)\right) \label{Qterm2}\\
& +  \frac{1}{2}\left(a_{max} + S_{\max} (|\mI|+1)\right)^2. \notag
\end{align}
In the final inequality, we use the fact that the term, \eqref{Qterm1}, is bounded by the term \eqref{Qterm2}.

Secondly, if $Q_i(t) \geq S_i(t)$ then, according to \eqref{Qcases},
\begin{align*}
 \frac{1}{2} Q_i(t+1)^2- \frac{1}{2} Q_i(t)^2 \leq & - Q_i(t) \left(   S_i(t) - A_i(t) - \sum_{i' : i'i\in\mL }S_{i'}(t)  p_{i'i}(t) \right) \\
& +\frac{1}{2} \left( S_i(t) - A_i(t) - \sum_{i' : i'i\in\mL }S_{i'}(t)  p_{i'i}(t) \right)^2\\
\leq & - Q_i(t) \left(   S_i(t) - A_i(t) - \sum_{i' : i'i\in\mL }S_{i'}(t)  p_{i'i}(t) \right) \\
& +\frac{1}{2} \left( a_{\max} + (|\mI|+1)S_{\max} \right)^2 .
\end{align*}
Thus defining
\begin{equation*}
K_0= \frac{1}{2} \left( a_{\max} + (|\mI|+1)S_{\max} \right)^2 + S_{\max} \left(a_{max} + S_{\max} (|\mI|+1)\right),
\end{equation*}
we see that in both cases, above, we have the required bound
\begin{equation*}
\frac{1}{2} Q_i(t+1)^2 - \frac{1}{2} Q_i(t)^2  \leq Q_i(t) \left( A_i(t) - S_i(t)
 +\sum_{i':i'i\in \mL} S_{i'}(t)p_{i'i}(t) \right)   + K_0.
\end{equation*}
\end{proof}

\begin{lemma}\label{QqqLemma}
There exists a constant $K_1>0$ such that the following equality holds
\begin{equation} \label{Qq}
\bE \left[   Q_i(t) S_{i'}(t)p_{i'i}(t)  \Big| Q(t),\delta(t)  \right]
\hspace{-4pt}\leq \bE \Big[ Q_i(t)\hspace{-8pt} \sum_{\sigma \in \mS_j(i')}\hspace{-8pt}
    \sigma_{i\rq{}}(t) P^{j(i')}_\sigma\bar{q}_{i'i}(t)  \Big| Q(t),\delta(t) \Big]+ K_1.
\end{equation}
\end{lemma}
\begin{proof}
First let us suppose that the queue has been empty over the last $k$ time steps. Then, since a bounded number of cars arrive at the queue per traffic cycle, the queue size $Q_i(t)$ must be less than $ (a_{\max} +  \sigma_{max} |\mL|)$. Clearly the above bound holds for any, $Q_i(t) \leq K_1= S_{\max}(a_{\max} +  \sigma_{max} |\mL|)$.
Now lets suppose $Q_i(t) \geq K_1$, we can take $Q_i(t)$ out the conditional expectation because it is known
\begin{equation*}
\bE[Q_i(t) S_{i'}(t) p_{i'i}(t)|Q(t),\delta(t)] = Q_i(t) \bE[ S_{i'}(t) p_{i'i}(t)|Q(t),\delta(t)].
\end{equation*}
The proportion of traffic $p(t)$ is independent of $S(t)$ and of $Q(t)$. So the
expectation of $p_{i'i}(t)$ given $Q(t)$ (and $S(t)$) is its mean $\bar{p}_{ii'}$. So
\begin{equation*}
\bE[ S_{i'}(t) p_{i'i}(t)|Q(t),\delta(t)] = \bE[ S_{i'}(t) \bar{p}_{i'i} |Q(t),\delta(t)] =
\bar{p}_{i'i} \bE[ S_{i'}(t)  |Q(t),\delta(t)].
\end{equation*}
Also \eqref{Scond} implies
\begin{equation*}
\bE [ S_{i'}(t) | Q(t),\delta(t) ] =\sum_{\sigma \in \mS_{j(i')}}\sigma_{i'} P^{j(i')}_\sigma.
\end{equation*}
Also since $q_{i'i}(t)$ is an unbiased estimate of $p_{i'i}$ at time $t$,
and independent of $Q(t)$ by assumption,
\begin{equation*}
\bar{p}_{i'i} = \bE[ q_{i'i}(t) |Q(t)].
\end{equation*}
Substituting this all back in, we have
\begin{align*}
\bE[Q_i(t) S_{i'}(t) p_{i'i}(t)|Q(t),\delta(t)]& = Q_i(t) \bar{p}_{i'i} \bE[ S_{i'}(t)|Q(t),\delta(t)] \\
& = Q_i(t) \bE[ q_{i'i}(t) |Q(t),\delta(t)] \sum_{\sigma \in \mS_{j(i')}} \sigma_{i'} P^{j(i')}_\sigma\\
& = \bE \Big[  Q_i(t) \sum_{\sigma \in \mS_j(i')} {\sigma}_{i\rq{}} P^{j(i')}_\sigma\bar{q}_{i'i}(t)  \Big| Q(t),\delta(t) \Big],
\end{align*}
thus the above inequality also holds in the case $Q_i(t)>K_1$  as required.
\end{proof}

Lemma \ref{swapsum} indicates an allowed reordering of terms.
\begin{lemma}\label{swapsum} The following equality holds for each measured queue size vector,
\begin{align*}
&  \sum_{i\in\mI}  \hat{Q}_i(t)  \left( \sum_{\sigma \in \mS_{j(i)}}  \sigma_i P^{j(i)}_{\sigma}(t)    -  \sum_{i':(i',i)\in \mL} \sum_{\sigma \in \mS_{j'(i')}} \sigma_{i'} P^{j(i')}_\sigma\bar{q}_{i'i}(t)\right) \\
&=   \sum_{k\in\mJ}\sum_{\sigma\in\mS_{k}}  P^{k}_{\sigma}(t)  \sum_{i\in k}  \sigma_i \left( \hat{Q}_i(t) -   \sum_{i': ii'\in\mL} \hat{Q}_{i'}(t) \bar{q}_{ii'}(t)  \right).
\end{align*}
\end{lemma}
\begin{proof}
Although the following set of equalities is some what lengthy, the premise is fairly simple. We want to change to order of summation so that we first sum over junctions $\mJ$ instead of first summing over in-roads $\mI$. These manipulations are as follows
\begin{align*}
&\sum_{i\in\mI} \hat{Q}_i(t)  \left( \sum_{\sigma \in \mS_{j(i)}} \sigma_i P^{j(i)}_{\sigma}(t)    -  \sum_{i':(i',i)\in \mL} \sum_{\sigma \in \mS_{j'(i')}}\sigma_{i'} P^{j(i')}_\sigma\bar{q}_{i'i}(t)\right) \\
=& \sum_{i\in\mI}  \sum_{k\in\mJ} \sum_{\sigma\in\mS_{k}} \hat{Q}_i(t) \sigma_i P^{k}_{\sigma}(t) \bI [ i\in k ] -\hspace{-3pt} \sum_{i\in\mI}  \sum_{i': i'i\in\mL} \sum_{k\in\mJ } \sum_{\sigma\in\mS_{k}}\hspace{-3pt} \hat{Q}_i(t)\sigma_{i'} P^{k}_\sigma \bar{q}_{i'i}(t) \bI [i' \hspace{-3pt}\in\hspace{-3pt} k] \\
=&  \sum_{k\in\mJ}\sum_{\sigma\in\mS_{k}}  P^{k}_{\sigma}(t) \hspace{-3pt}  \sum_{i\in\mI} \hat{Q}_i(t) \sigma_i \bI [i \in k ]
- \hspace{-3pt}  \sum_{k\in\mJ } \sum_{\sigma\in\mS_{k}}  P^{k}_\sigma \sum_{i\in\mI}  \sum_{i': i'i\in\mL}\hspace{-3pt} \hat{Q}_i(t) \sigma_{i'} \bar{q}_{i'i}(t) \bI [i'\hspace{-3pt} \in\hspace{-3pt} k]\\
=&  \sum_{k\in\mJ}\sum_{\sigma\in\mS_{k}}  P^{k}_{\sigma}(t)   \left( \sum_{i\in k} \hat{Q}_i(t) \sigma_i -  \sum_{i'\in\mI}  \sum_{i: i'i\in\mL} \hat{Q}_i(t) \sigma_{i'} \bar{q}_{i'i}(t) \bI [i' \in k ]  \right)\\
=&  \sum_{k\in\mJ}\sum_{\sigma\in\mS_{k}}  P^{k}_{\sigma}(t)  \sum_{i\in k}  \sigma_i  \left( \hat{Q}_i(t) -   \sum_{i': ii'\in\mL} \hat{Q}_{i'}(t) \bar{q}_{ii'}(t)  \right).
\end{align*}
In the first equality above, we expand brackets. In the second equality, we reorder the summation so we first sum over junctions and then over schedules. In the third and fourth equality, we collect together terms for each in-road.
\end{proof}

Our last lemma gives a bound on an optimization problem to be used later.
\begin{lemma}\label{epslemma} If ${a} + \epsilon\mathbf{1}\in \mA$  then
\begin{equation*}
  \epsilon
  < \left(1-\frac{L}{T}\right) \min_{u\geq 0: \sum_{i\in \mI}u_i =1}
    \left(
        \sum_{j\in\mJ}\max_{\sigma\in\mS_j}
            \left\{
            w_{\sigma}\left( u \right)
                \right\}
        -
        \sum_{i\in\mI}   u_i\bar{a}_i(t)
    \right) .
\end{equation*}
\end{lemma}
\begin{proof}
 By definition ${a} + \epsilon\mathbf{1}\in \mA$ when
 there exists a positive vector $\rho=(\rho^j_{\sigma}:  \sigma\in\mS_j, j\in\mJ)$ and a positive vector $s=(s_i:i\in\mI)$ satisfying the constraints
\begin{align}
&a_i + \epsilon  + \sum_{i' : i'i\in\mL } s_{i'}\bar{p}_{i'i}   < s_i , \label{aequn} \\
\sum_{\sigma\in\mS_j} & \rho^j_{\sigma}   \leq 1-\frac{L}{T},\qquad s_i \leq \sum_{\sigma\in\mS_j} \rho^j_{\sigma}\sigma_i \label{aequn2}
\end{align}
\text{for each} $j\in\mJ$\text{ and }$i\in j$. We can express \eqref{aequn} more concisely in vector form as $a + \epsilon \mathbf{1} < s^{\textsf{T}} (I - \bar{p})$.
Notice the inverse of $(I-\bar{p})$ is the positive matrix $(I - \bar{p})^{-1} = I + \bar{p} + \bar{p}^2 + ...$ Thus we can equivalently express condition \eqref{aequn}  as
\begin{equation}
(a + \epsilon \mathbf{1})^{\textsf{T}} (I - \bar{p})  < s^{\textsf{T}}.
\end{equation}
We can now observe that if we replace $s$ as above with
\begin{equation}
\tilde{s}_i = \sum_{\sigma\in\mS_j} \rho^j_{\sigma}\sigma_i
\end{equation}
then equations (\ref{aequn}-\ref{aequn2}) must hold. In other words there exists a  $\rho=(\rho^j_{\sigma}:  \sigma\in\mS_j, j\in\mJ)$ such that
\begin{align}
a_i +\epsilon + \sum_{i' : i'i\in\mL } \sum_{\sigma\in\mS_{j(i')}} \rho^{j(i')}_{\sigma} \sigma_{i'} \bar{p}_{i'i}   & < \sum_{\sigma\in\mS_j} \rho^j_{\sigma}\sigma_i  ,&& \text{for } i\in j,\text{ with } j\in\mJ,\label{Lemma:Stab1} \\
\sum_{\sigma\in\mS_j} & \rho^j_{\sigma}   \leq 1-\frac{L}{T}, && \text{for } j\in\mJ. \label{Lemma:Stab2}
\end{align}
We now focus on the inequality \eqref{Lemma:Stab1}. Rearranging it, the above holds when there exists $\rho$ such that
\begin{align}
\epsilon   & < \min_{j\in\mJ, i\in j}\left\{  \sum_{\sigma\in\mS_j} \rho^j_{\sigma}\sigma_i - \sum_{i' : i'i\in\mL } \sum_{\sigma\in\mS_{j(i')}} \rho^{j(i')}_{\sigma} \sigma_{i'} \bar{p}_{i'i}  - a_i  \right\} ,&& \label{Lemma:Stab3} \\
& \qquad\qquad\qquad \sum_{\sigma\in\mS_j}  \rho^j_{\sigma}   \leq 1-\frac{L}{T}, && j\in\mJ.\label{Lemma:Stab4}
\end{align}
The above statement can only hold if it also holds when we maximize over $\rho$ thus the following must hold
\begin{align} \label{Lemma:Stab5}
\epsilon    < \max_{\rho:\forall j\in\mJ,\sum_{\sigma\in\mS_j} \rho^j_{\sigma}   \leq 1-L/T} \min_{j\in\mJ,i\in j}\left\{  \sum_{\sigma\in\mS_j} \rho^j_{\sigma}\sigma_i - \sum_{i' : i'i\in\mL } \sum_{\sigma\in\mS_{j(i')}} \rho^{j(i')}_{\sigma} \sigma_{i'} \bar{p}_{i'i}  - a_i  \right\}.
\end{align}
Since the minimum of a finite set in $\bR$ is equal to
the minimum of the convex combinations of elements in that set,
\begin{align}
 &\min_{i\in \mI}\left\{  \sum_{\sigma\in\mS_{j(i)}} \rho^{j(i)}_{\sigma}\sigma_i - \sum_{i' : i'i\in\mL } \sum_{\sigma\in\mS_{j(i')}} \rho^{j(i')}_{\sigma} \sigma_{i'} \bar{p}_{i'i}  - a_i  \right\} \notag \\
= &\min_{u\geq 0: \sum_{i\in \mI}u_i =1} \left\{ \sum_{i\in \mI}
\sum_{\sigma\in\mS_{j(i)}} \rho^{j(i)}_{\sigma}\sigma_i u_i
- \sum_{i\in \mI} \sum_{i' : i'i\in\mL } \sum_{\sigma\in\mS_{j(i')}} \rho^{j(i')}_{\sigma} \sigma_{i'} \bar{p}_{i'i}u_i
 - \sum_{i\in \in\mI} a_iu_i  \right\}. \label{minEqul}
\end{align}
Next by exactly the same argument used to prove Lemma \ref{swapsum}, we have that
\begin{align*}
  &\sum_{i\in \mI} \sum_{\sigma\in\mS_{j(i)}} \rho^{j(i)}_{\sigma}\sigma_i u_i
  - \sum_{i\in \mI}  \sum_{i' : i'i\in\mL } \sum_{\sigma\in\mS_{j(i')}} \rho^{j(i')}_{\sigma} \sigma_{i'} \bar{p}_{i'i}u_i  -   \sum_{i\in \mI}  a_iu_i \\
=& \sum_{k\in\mJ}\sum_{\sigma\in\mS_{k}}  \rho^{k}_{\sigma}  \sum_{i\in k}  \sigma_i  \left(   u_i -   \sum_{i': ii'\in\mL} u_{i'} \bar{p}_{ii'} \right) -   \sum_{i\in \mI}  a_iu_i.
\end{align*}
Substituting this equality into \eqref{minEqul}, we have that \eqref{Lemma:Stab5} reads as
\begin{align}
\epsilon    < \max_{\rho:\forall j\in\mJ,\sum_{\sigma\in\mS_j} \rho^j_{\sigma}   \leq 1-L/T } \min_{u\geq 0: \sum_{i\in \mI}u_i =1}\left\{ \sum_{k\in\mJ}\sum_{\sigma\in\mS_{k}}  \rho^{k}_{\sigma}  \sum_{i\in k}  \sigma_i  \left(   u_i -   \sum_{i': ii'\in\mL} u_{i'} \bar{p}_{ii'} \right) -   \sum_{i\in \mI}  a_iu_i \right\}.
\end{align}
Finally for any function $f(u,\rho)$, it holds that $\min_u \max_\rho f(u,\rho) \geq \max_\rho \min_u  f(u,\rho)$. Thus if $\bar{a} + \epsilon\mathbf{1}\in \mA$ then it must be true that
\begin{equation*}
\epsilon    < \min_{u\geq 0: \sum_{i\in \mI}u_i =1}  \max_{\rho : \forall
j\in\mJ,\sum_{\sigma\in\mS_j} \rho^j_{\sigma}   \leq 1-L/T} \left\{    \sum_{k\in\mJ}\sum_{\sigma\in\mS_{k}}  \rho^{k}_{\sigma}  \sum_{i\in k}  \sigma_i  \left(   u_i -   \sum_{i': ii'\in\mL} u_{i'} \bar{p}_{ii'} \right) -   \sum_{i\in \mI}  a_iu_i  \right\}.
\end{equation*}
Finally, we note that the above maximization over $\rho$ must be achieved at a value where, for each $j\in\mJ$, $\rho^j_\sigma=1-\frac{L}{T}$ for some $\sigma$, so
\begin{equation*}
    \epsilon
    <\left(1-\frac{L}{T}\right) \min_{u\geq 0: \sum_{i\in \mI}u_i =1}  \left\{
        \sum_{j\in\mJ} \max_{\sigma\in\mS_j}\sum_{i\in j}
		\sigma_i  \left(
			u_i -   \sum_{i': ii'\in\mL} u_{i'} \bar{p}_{ii'}
		    \right)
	- \sum_{i\in \mI}  a_iu_i
    \right\}.
\end{equation*}

\end{proof}

\subsection{Proofs}\label{sec:apx.proof}
We now provide a proof that the road traffic network is unstable whenever the arrival rates are outside the set $\bar{\mA}$.

\begin{proof}[Proof of Proposition \ref{prop1}]
If $\bar{a}\notin \bar{\mA}$ then there exists an $\epsilon>0$ where, for any vector $\rho=(\rho_\sigma^j : \sigma \in\mS_j, j\in\mJ)$ and $s=(s_i:i\in\mI)$, satisfying
\begin{align}\label{rhosum}
\sum_{\sigma\in\mS_j} & \rho^j_{\sigma}   \leq 1-\frac{L}{T},\qquad\text{and}\qquad s_{i'} \leq \sum_{\sigma\in\mS_j} \rho^j_{\sigma}\sigma_{i'}
\end{align}
for $ j\in\mJ$ and ${i'}\in\mI$, and there exists an in-road $i$ such that
\begin{align}\label{epsibound}
a_i  -s_i  + \sum_{i' : i'i\in\mL }s_{i'} \bar{p}_{i'i}   & > \epsilon .
\end{align}

Consider any policy $P_\sigma$, and consider the average service devoted to each queue and the number of departures from each queue
\begin{equation*}
\rho_\sigma^j(\tau) = \frac{1}{\tau} \sum_{t=1}^{\tau} P_{\sigma}^j(t)\quad \text{and} \quad s_i(\tau)=\frac{1}{\tau} \sum_{t=1}^{\tau} S_i(t).
\end{equation*}
Taking a suitable subsequence if necessary, the sequences $\rho_\sigma^j(t)$ and $s_i(t)$  must converge to some value $\rho_\sigma^j$ and $s_i$ satisfying condition \eqref{rhosum}. The long run queue size must converge to the average arrival minus the average departures. Thus we have that
\begin{equation*}
\lim_{t\rightarrow\infty} \frac{Q^\Sigma(t)}{t} = \sum_{i\in\mI} \left( \bar{a}_i  + \sum_{i' : i'i\in\mL } s_{i'} \bar{p}_{i'i} - s_i \right) > \epsilon.
\end{equation*}
In the above inequality, we note that each term in the above summation is positive since it is the limit of positive queue sizes and one of those terms in greater than $\epsilon$ by \eqref{epsibound}.
Thus we see that for any policy there always exists a queue that is unstable: there exists a $t_0$ such that for all $t> t_0$
\begin{equation*}
Q^\Sigma(t) > t \epsilon,
\end{equation*}
thus,
\begin{equation*}
\lim_{\tau\rightarrow \infty} \frac{1}{\tau} \sum_{t=1}^{\tau} Q^\Sigma(t) \geq \lim_{\tau\rightarrow \infty} \frac{1}{\tau} \sum_{t=t_0}^{\tau} \epsilon t = \infty.
\end{equation*}
So, taking expectations, we see that the network must be unstable:
\begin{equation*}
\limsup_{\tau\rightarrow \infty} \bE \left[ \frac{1}{\tau}  \sum_{t=1}^{\tau} Q^\Sigma(t) \right] \geq  \bE \left[ \liminf_{\tau\rightarrow\infty} \frac{1}{\tau}  \sum_{t=1}^{\tau} Q^\Sigma(t) \right]  = \infty.
\end{equation*}
In the last inequality we apply Fatou's Lemma.
\end{proof}

\subsection{Formal proof of Theorem \ref{thrm1}}
We now begin to develop the proof of Theorem \ref{thrm1}. First, we require a bound on the change in the (euclidean) distance of our queue sizes from zero. This is proven in the following proposition.

\begin{proposition} \label{QProp} There exists a constant \hl{$K^*$ }such that
\begin{align*}
& \sum_{i\in\mI} \bE\left[ \frac{1}{2} Q_i(t+1)^2- \frac{1}{2} Q_i(t)^2 \bigg| Q(t),\delta(t) \right]\\
\leq &\bE\left[ \sum_{i\in\mI}   Q_i(t) \bar{a}_i(t)  -  \sum_{j\in\mJ}\sum_{\sigma \in \mS_{j}}    P^{j}_\sigma   w_{\sigma}(Q(t)) \Bigg| Q(t),\delta(t) \right] +K^*
\end{align*}
where $w_{\sigma}(Q(t))$ is defined by \eqref{wdef}.
\end{proposition}

\begin{proof}
\noindent We can expand the left side through the following inequalities to reach the desired bound
\begin{align*}
&\sum_{i\in\mI} \bE\left[  \frac{1}{2} Q_i(t+1)^2 - \frac{1}{2} Q_i(t)^2  \bigg| Q(t),\delta(t) \right] \notag \\
  \leq  &\sum_{i\in\mI} \bE\left[  Q_i(t) \left( A_i(t) - S_i(t)
 +\sum_{i':i'i\in \mL} S_{i'}(t)p_{i'i}(t) \right)  \bigg| Q(t),\delta(t) \right]  + K_0 \\
 \leq &  \sum_{i\in\mI}  \bE\hspace{-3pt}\left[ Q_i(t)\hspace{-5pt}
    \left(
        \bar{a}_i(t)
        -\hspace{-10pt} \sum_{\sigma \in \mS_{j(i)}}\hspace{-6pt} {\sigma}_{i}  P^{j(i)}_{\sigma}(t)
        +\hspace{-10pt}\sum_{i':i'i\in \mL} \sum_{\sigma \in \mS_{j(i')}}
            \hspace{-6pt}{\sigma}_{i'} P^{j(i')}_\sigma\bar{q}_{i'i}(t)
    \right)\hspace{-6pt}
    \Bigg| Q(t),\delta(t) \right] \hspace{-3pt}+\hspace{-3pt} \tilde{K}\notag\\
    = & \bE\left[ \sum_{i\in\mI}   Q_i(t) \bar{a}_i(t) \Bigg| Q(t),\delta(t) \right] +\tilde{K}
                    \notag\\
    - & \bE\hspace{-3pt}\left[ \hspace{-3pt} \sum_{i\in\mI}\hspace{-3pt} \left(\hat{Q}_i(t)-\delta_i(t)\right)\hspace{-6pt}  \left(
        \hspace{-3pt}\sum_{\sigma \in \mS_{j(i)}}\hspace{-6pt} {\sigma}_{i} P^{j(i)}_{\sigma}(t)
        - \hspace{-6pt} \sum_{i':i'i\in \mL} \sum_{\sigma \in \mS_{j(i')}}\hspace{-6pt}
            {\sigma}_{i'}  P^{j(i')}_\sigma\bar{q}_{i'i}(t)
        \right)\hspace{-3pt} \Bigg| Q(t),\delta(t) \right]\\
\leq &\bE\left[ \sum_{i\in\mI}   Q_i(t) \bar{a}_i(t)   -  \sum_{j\in\mJ}\sum_{\sigma \in \mS_{j}}    P^{j}_\sigma   w_{\sigma }(\hat{Q}(t)) \Bigg| Q(t),\delta(t) \right]+K^*.
\end{align*}

The first inequality can be reached by using Lemma \ref{QKLemma} to expand the recursion \eqref{Qprocess} for some constant $K_0>0$. Then by using the inequalities
\begin{equation}  \label{Qqq}
\bE\hspace{-3pt} \left[   Q_i(t) S_{i'}(t)p_{i'i}(t)  \Big| Q(t),\delta(t)  \right]
\hspace{-3pt} \leq\hspace{-3
pt} \bE \Big[  Q_i(t) \sum_{\sigma \in \mS_{j(i')}} {\sigma}_{i'} P^{j(i')}_\sigma\bar{q}_{i'i}(t)  \Big| Q(t),\delta(t) \Big]\hspace{-3pt} + K_1
\end{equation}
and
\begin{equation}
\bE\left[ Q_i(t)  A_i(t) \Big| Q(t) \right] = Q_i(t) \bar{a}_i(t) \label{Qaa}
\end{equation}
with rearranging we can further expand by taking constant $\tilde{K}=K_0 + |\mL| K_1$ for instance. \eqref{Qqq} is proven in Lemma \ref{QqqLemma}, whereas \eqref{Qaa} holds by definition \eqref{abar}.
The last inequality is given by swapping the order of summation inside the expectation from a summation over queues on in-roads $\mI$ and then as a summation over schedules $\mS_j$ to the other way around, and using the bound on the error terms in \eqref{error_bound}.
By Lemma \ref{swapsum} and by the definition of the weights $w_\sigma(\hat{Q}(t))$ in \eqref{wdef}, we have that
\begin{align*}
&  \sum_{i\in\mI} \hat{Q}_i(t)  \left( \sum_{\sigma \in \mS_{j(i)}}  {\sigma}_{i}
P^{j(i)}_{\sigma}(t)    -  \sum_{i':i'i\in \mL} \sum_{\sigma \in \mS_{j(i')}} {\sigma}_{i'} P^{j(i')}_\sigma\bar{q}_{i'i}(t)\right) \\
= &   \sum_{k\in\mJ}\sum_{\sigma\in\mS_{k}}  P^{k}_{\sigma}(t)  \sum_{i\in k} {\sigma}_{i}  \left( \hat{Q}_i(t) -   \sum_{i': ii'\in\mL} \hat{Q}_{i'}(t) \bar{q}_{ii'}(t)  \right)\\
= & \sum_{k\in\mJ} \sum_{\sigma \in \mS_{k}}    P^{k}_\sigma   w_{\sigma }(\hat{Q}(t)).
\end{align*}
\hl{Furthermore we can bound the extra term arising from the measurement error as
\begin{align*}
&\bE\left[  \sum_{i\in\mI} \delta_i(t) \left(\sum_{\sigma \in \mS_{j(i)}} {\sigma}_{i} P^{j(i)}_{\sigma}(t)-  \sum_{i':i'i\in \mL} \sum_{\sigma \in \mS_{j(i')}}{\sigma}_{i'}  P^{j(i')}_\sigma\bar{q}_{i'i}(t)\right) \Bigg| Q(t),\delta(t) \right]\\
=& \sum_{i\in\mI} \delta_i(t)\bE\left[\left(\sum_{\sigma \in \mS_{j(i)}} {\sigma}_{i} P^{j(i)}_{\sigma}(t)-  \sum_{i':i'i\in \mL} \sum_{\sigma \in \mS_{j(i')}}{\sigma}_{i'}  P^{j(i')}_\sigma\bar{q}_{i'i}(t)\right) \Bigg| Q(t),\delta(t) \right]\\
\leq&\delta_{\max}\cdot |\mI|\cdot S_{\max} (1+|\mL|),
\end{align*}
which by introducing $K^*=\tilde{K}+\delta_{\max}\cdot |\mI|\cdot S_{\max} (1+|\mL|)$ completes our proof.}
\end{proof}

Now that we have proven the previous proposition, we are able to prove the main mathematical result of this paper, Theorem \ref{thrm1}.

\begin{proof}[Proof of Theorem \ref{thrm1}]
By Proposition \ref{QProp}, we have that
\begin{align}\label{Lya2}
&\sum_{i\in\mI} \bE\bigg[ \frac{1}{2} Q_i(t+1)^2- \frac{1}{2} Q_i(t)^2 \bigg| Q(t),\delta(t) \bigg] \notag \\
\leq &\bE\bigg[ \sum_{i\in\mI}   Q_i(t) \bar{a}_i(t)   -  \sum_{j\in\mJ}\sum_{\sigma \in \mS_{j}}    P^{j}_\sigma   w_{\sigma }(\hat{Q}(t))\bigg| Q(t),\delta(t) \bigg] +K^*.
\end{align}
Further, by Lemma \ref{weightslemma2}, we know that
\begin{equation}
 \sum_{\sigma \in \mS_{j}}    P^{j}_\sigma   w_{\sigma }(\hat{Q}(t)) \geq  \max_{\sigma\in\mS_j}\Big\{   w_{\sigma}(Q(t)) \Big\}-\delta_{\max}S_{\max}|\mI| - \frac{1}{\eta} \log |\mS_j|.
\end{equation}
Applying this bound to \eqref{Lya2} and taking expectations, we have that
\begin{equation}
    \sum_{i\in\mI} \bE\hspace{-3pt}\left[
        \frac{1}{2} Q_i(t+1)^2- \frac{1}{2} Q_i(t)^2\right]
\hspace{-3pt}\leq\hspace{-3pt}  \bE\hspace{-3pt}\left[ \sum_{i\in\mI}   Q_i(t) \bar{a}_i(t)
   -\hspace{-3pt} \sum_{j\in\mJ}  \max_{\sigma\in\mS_j}  \left\{
    w_{\sigma}(Q(t))
        \right\}  \right]\hspace{-3pt}+ K
\label{QWbounder}
\end{equation}
\hl{where
\begin{equation*}
K = K^*+\delta_{\max}S_{\max}|\mI| + \sum_{j\in\mJ} \frac{1}{\eta}\log |\mS_j|.
\end{equation*}}
We focus on bounding the term \eqref{QWbounder} above. Recalling the definition of $Q^\Sigma(t)$, \eqref{QSigma},
observe that
\begin{align}
    &\bE \left[  \sum_{j\in\mJ}\max_{\sigma\in\mS_j}  \left\{
    w_{\sigma}(Q(t))
        \right\} - \sum_{i\in\mI}   Q_i(t) \bar{a}_i(t) \right]    \notag\\
= & \bE \left[  Q^\Sigma(t)\times
 \left(  \sum_{j\in\mJ}\max_{\sigma\in\mS_j}  \left\{
    w_{\sigma}\left( \frac{Q(t)}{Q^\Sigma (t)} \right)
        \right\} - \sum_{i\in\mI}   \frac{Q_i(t)}{Q^\Sigma (t)} \bar{a}_i(t) \right) \right]  \notag\\
\geq &  \bE \left[  Q^\Sigma(t) \times \min_{u\geq 0: \sum_{i\in \mI}u_i =1} \left(  \sum_{j\in\mJ}\max_{\sigma\in\mS_j}  \left\{
    w_{\sigma}\left( u \right)
        \right\} - \sum_{i\in\mI}   u_i\bar{a}_i(t) \right)  \right]\notag\\
\geq &  \bE \left[  Q^\Sigma(t)  \times \epsilon\right]. \label{epsQsig}
\end{align}
In the first inequality above, we substitute $u_i$ for
$\frac{Q_{i}(t)}{Q^\Sigma(t)}$
and minimize over $u_i$. In the second inequality, we apply our bound from Lemma \ref{epslemma}:
\hl{\begin{equation*}
  \epsilon
  <\left(1+\frac{L}{T}\right)  \min_{u\geq 0: \sum_{i\in \mI}u_i =1}
    \left(
        \sum_{j\in\mJ}\max_{\sigma\in\mS_j}
            \left\{
            w_{\sigma}\left( u \right)
                \right\}
        -
        \sum_{i\in\mI}   u_i\bar{a}_i(t)
    \right) .
\end{equation*}}
Applying \eqref{epsQsig} to \eqref{QWbounder}, we gain a far simpler bound
\begin{align}
&\sum_{i\in\mI} \bE\left[ \frac{1}{2} Q_i(t+1)^2- \frac{1}{2} Q_i(t)^2  \right] \leq - \epsilon\bE \left[  Q^\Sigma(t)  \right]  +K.
\end{align}
Summing over $t=0,...,\tau-1$, we have
\begin{equation}
\sum_{i\in\mI} \bE\left[ \frac{1}{2} Q_i(\tau)^2- \frac{1}{2} Q_i(0)^2  \right]
\leq -\epsilon \bE \left[ \sum_{t=0}^{\tau-1}  Q^\Sigma(t)  \right]  + \tau K.
\end{equation}
Finally, rearranging and dividing by \hl{$\tau$} gains the required bound,
\begin{align*}
\lim_{\tau\rightarrow\infty}
    \bE \left[ \frac{1}{\tau} \sum_{t=0}^{\tau-1}  Q^\Sigma(t)  \right]
&\leq
    \frac{K}{\epsilon}
    +
    \lim_{\tau\rightarrow\infty} \frac{1}{\tau\epsilon}
            \bE\left[ \frac{1}{2} \sum_{i\in\mI} \  Q_i(0)^2 \right] = \frac{K}{\epsilon}  < \infty.
\end{align*}

\end{proof}

\bibliographystyle{elsarticle/model2-names}

\bibliography{ref}

\end{document}